\documentclass[12pt]{article}
\usepackage{amssymb,amsmath}
\usepackage{color}

\newcommand{\R}{\mathbb{R}}
\newcommand{\N}{\mathbb{N}}
\newcommand{\F}{\mathbb{F}}
\def\half{{1\over 2}}
\def\abs#1{|#1|}

\def\norm#1{\|#1\|}
\def\cat{\mathop{\rm cat}\nolimits}
\def\dist{\mathop{\rm dist}\,}

\def\claim#1.{\noindent {\bf #1.}}

\def\flushright#1{{\unskip\nobreak\hfil\penalty50\hskip2em\hbox{}\nobreak\hfil%
#1\parfillskip=0pt\finalhyphendemerits=0\par}}

\def\bull{\vrule height 1.8ex width 1.0ex depth .1ex }
\def\QED{\ifmmode\eqno\hbox{$\bull$}\else\flushright{\hbox{$\bull$}}\fi}

\def\wlimit{\rightharpoonup}

\def\supp{\mathop{\rm supp}}

\def\be{\begin{equation}}
\def\ee{\end{equation}}

\def\bea{\begin{eqnarray*}}
\def\eea{\end{eqnarray*}}

\def\epsilon{\varepsilon}
\def\intRN{\int_{\R^N}}

\def\Je{J_\epsilon}

\def\calV{{\cal V}}
\def\calK{{\cal K}}
\def\calD{{\cal D}}

\def\whS{\widehat S}

\def\wrho{\widehat\rho}

\def\cuplength{{\rm cupl}}

\def\hdelta{\hat\delta}
\def\wJ{\widetilde J}

\def\calXed{{\cal X}_{\epsilon}}
\def\calXp{\calXed^{E(m_0)+\hdelta}}
\def\calXm{\calXed^{E(m_0)-\hdelta}}

\def\Bonec{B(0,{1\over\sqrt\e})^c}
\def\Btwoc{B(0,{2\over\sqrt\e})^c}

\def\whY{\widehat Y_\e}

\newcommand{\e}{\varepsilon}

\def\zR{\widetilde \zeta_{{1\over\sqrt\e}+{k\over\e^{1/4}}}}

\newtheorem{theorem}{Theorem}[section]
\newtheorem{proposition}[theorem]{Proposition}
\newtheorem{corollary}[theorem]{Corollary}
\newtheorem{lemma}[theorem]{Lemma}
\newtheorem{remark}[theorem]{Remark}
\newtheorem{definition}[theorem]{Definition}

\begin{document}

\title 
{Semi-classical states for the nonlinear Choquard equations: 
existence, multiplicity and concentration at a potential well}

\author{Silvia Cingolani\\
Dipartimento di  Matematica\\
Universit\`{a} degli Studi di Bari Aldo Moro
\\
Via E. Orabona 4, 
70125 Bari, Italy\\
\\
Kazunaga Tanaka\\
Department of Mathematics\\
School of Science and Engineering\\
Waseda University\\
3-4-1 Ohkubo, Shijuku-ku, Tokyo 169-8555, Japan}

\date{}

\maketitle


\begin{abstract}
We study existence and multiplicity of semi-classical states for the nonlinear Choquard equation:
    $$  -\e^2\Delta v+V(x)v = {1\over\e^\alpha}(I_\alpha*F(v))f(v) \quad \hbox{in}\ \R^N,
    $$
where $N\geq 3$, $\alpha\in (0,N)$, $I_\alpha(x)={A_\alpha\over |x|^{N-\alpha}}$ is the Riesz potential,
$F\in C^1(\R,\R)$, $F'(s) = f(s)$ and $\e>0$ is a small parameter.  

We develop a new variational approach and we show the existence of a family of solutions concentrating,
as $\e\to 0$, to a local minima of $V(x)$ under general conditions on $F(s)$.
Our result is new also for $f(s)=|s|^{p-2}s$ and applicable for $p\in ({N+\alpha\over N}, {N+\alpha\over N-2})$.
Especially, we can give the existence result for locally sublinear case $p\in ({N+\alpha\over N},2)$, which
gives a positive answer to an open problem arisen in recent works of Moroz and Van Schaftingen.

We also study the multiplicity of positive single-peak solutions and we show the existence of 
at least $\cuplength(K)+1$ solutions concentrating around $K$ as $\e\to 0$, where $K\subset \Omega$ is the
set of minima of $V(x)$ in a bounded potential well $\Omega$, that is, $m_0 \equiv \inf_{x\in \Omega} V(x) < \inf_{x\in \partial\Omega}V(x)$
and $K=\{x\in\Omega;\, V(x)=m_0\}$.
\end{abstract}

  \medskip




\section{\label{section:1} Introduction}
In this paper we consider the nonlocal equation
    \begin{equation} \label{eq:1.1}
    - \epsilon^2 \Delta v + V(x)v = \frac{1}{\epsilon^\alpha}  (I_\alpha * F(v)) f(v) \  \text{in $\R^N$},
    \quad v>0, \quad v \in H^1(\R^N),
    \end{equation}
where $\e>0$ is a positive parameter, $N\geq 3$, $\alpha\in (0,N)$, $I_\alpha:\, \R^N\setminus\{ 0\}\to \R$ the Riesz potential
defined by
    $$  I_\alpha(x)= \frac{\Gamma(\frac{N-\alpha}{2})}{\Gamma(\frac{\alpha}{2}) \pi^{N/2} 2^\alpha |x|^{N- \alpha}}
    $$
and $F(s)\in C^1(\R,\R)$ and $f(s)=F'(s)$.

In 1954 the equation \eqref{eq:1.1} with $N=3$, $\alpha=2$ and $F(s)=\half |s|^2$ was introduced by Pekar in \cite{pekar} to describe 
the quantum theory of a polaron at rest and in 1976 it was arisen in the work of Choquard on the modeling of an electron trapped in its own hole, 
in a certain approximation to Hartree -- Fock theory of one-component plasma  \cite{lieb}  (see also \cite{fl,fty}). 
In particular the equation \eqref{eq:1.1} with $V\equiv 1$ 
corresponds to the stationary nonlinear  Hartree equation.
Indeed if $v$ is a solution of 
\eqref{eq:1.1}, then the wave function $\psi(x,t) = e^{i t/\epsilon} v(x)$ is a solitary wave of the time-dependent Hartree equation  
\begin{equation} \label{eq:1.2}
\epsilon i \psi_t = - \epsilon^2 \Delta \psi  - \frac{1}{\epsilon^2}  ( \frac{1}{|x|} *  |\psi|^2)\psi \quad  \text{in $\R \times \R^3$}.
\end{equation}
Finally, we recall that the three-dimensional nonlocal equation  was also proposed by Penrose
\cite{pe1, pe2, pe3, mpt} in his discussion on the self-gravitational collapse of a quantum mechanical wave-function and in that context it is known as the Schr\"odinger-Newton equation.  Indeed, after a suitable change of variable, it can be written down as 
    \begin{equation}\label{eq:1.3}
        \begin{cases}
        -\e^2 \Delta \phi + V(x) \phi = U \phi , \\
        - \Delta U =  4 \pi \gamma |\phi|^2
        \end{cases}
    \end{equation}
where $\phi:\, \R^3\to \mathbb{C}$, $U:\, \R^3\to\R$ are unknown functions.

In literature the nonlocal equation \eqref{eq:1.1} is usually called 
the {\sl nonlinear Choquard equation}.

In the present paper we are interested to the search of semi-classical states for the nonlinear Choquard equation, namely
solutions to \eqref{eq:1.1} which exist for small values of $\epsilon$ and have a specific behavior as $\e$ tends to $0$.  We recall that the semiclassical analysis for nonlinear Schr\"odinger equations has been largely investigated in the last years, starting from the pioneering paper of Floer and Weistein \cite{FW} (see also \cite{ABC,DF, Oh1}).

Concerning the nonlocal case, Wei and Winter  \cite{WW}  established the existence of a family of solutions 
of the Schr\"odinger-Newton equations \eqref{eq:1.3}, concentrating as $\e\to 0$
around a non-degenerate critical point of the potential $V(x)$ when $\inf V(x)>0$ via a Lyapunov-Schmidt reduction method
(see also \cite{S} when $\liminf V(x) |x|^\gamma >0$ for some $\gamma \in [0,1)$).
We note that the Lyapunov-Schmidt reduction method relies on the existence and uniqueness (up to translation) of the positive
ground state solution of the limit problem:
\begin{equation} \label{eq:1.4}
  -\Delta u+u=({1\over|x|}*|u|^2)u \quad \hbox{in}\ \R^3,
\end{equation}
which is established by Lieb in \cite{lieb} (see also \cite{Lenzmann2009,l2}). 
Successively in \cite{ccs1,CSS} variational methods were employed to study the semiclassical analysis for the three-dimensional magnetic Choquard equation 
with $\alpha=2$ and $F(s)= \half |s|^2$.

In this work by means of a new approach we study the nonlinear Choquard equation \eqref{eq:1.1}, where $\alpha\in (0,N)$ and $f(s)$ is a general nonlinearity
which does not satisfy a monotonicity assumption nor Ambrosetti-Rabinowitz type condition \cite{AmbrosettiRabinowitz}.
Defining $u(x)=v(\e x)$, the equation \eqref{eq:1.1} is equivalent to
    \begin{equation}\label{eq:1.5}
    -\Delta u  +V(\e x) u =   (I_\alpha * F(u)) f(u) \ \text{in $\R^N$}, \  u > 0, \ u \in H^1(\R^N).
    \end{equation}
Thus, we try to find critical points of the corresponding functional:
    \begin{equation}\label{eq:1.6}
    J_\e(u)=\half \int_{\R^N} |\nabla u|^2+ V(\e x)u^2\, dx -\half \int_{\R^N} (I_\alpha*F(u))F(u)\, dx:\, H^1(\R^N)\to \R.
    \end{equation}
It is natural to ask the existence of a family $(u_\e)$ of solutions of \eqref{eq:1.5}, which concentrates to a local minima of
$V(x)$ as $\e$ tends to $0$.

To give our main result, we assume that

\vskip2pt

\begin{itemize}
\item[(f1)] $f \in C(\R, \R)$;
\item[(f2)] there exists $C >0$ such that for every $s \in \R$, $$|s f(s)| \leq C (|s|^{\frac{N + \alpha}{N}} + |s|^{\frac{N + \alpha}{N -2}});$$
\item[(f3)] $F(s)= \int_0^s f(t) dt$ satisfies
  $$\lim_{ s \to 0} \frac{F(s)}{|s|^{\frac{N+ \alpha}{N}}} =0, \quad
  \lim_{ s \to  + \infty} \frac{F(s)}{|s|^{\frac{N+ \alpha}{N-2}}} =0;$$
\item[(f4)] there exists $s_0 \in \R$, $s_0 \neq 0$ such that $F(s_0) > 0$.
\item[(f5)] $f$ is odd and $f$ is positive on $(0, +\infty)$.
\end{itemize}
\vskip 2pt
We remark that the conditions (f1)--(f5) are in the spirit of Berestycki and Lions (see \cite{BL,HIT,MS1}) and
cover the case $f(s)= |s|^{p-2} s$ with $p \in ({N+\alpha\over N}, {N+\alpha\over N-2})$.
\noindent
For $V(x)$, we assume
\begin{itemize}
\item[(V1)] $ V \in C(\R^N, \R)$, $\inf_{x \in \R^N} V(x) \equiv \underline{V} >0$, $\sup_{x \in \R^N} V(x) < \infty$.
  \noindent
\item[(V2)] There is a bounded domain $\Omega \subset \R^N$ such that
  \[ m_0 \equiv \inf_{x \in \Omega} V(x) <    \inf_{x \in \partial \Omega} V(x). \]
  We also set
  \begin{equation}\label{eq:1.7}
    K=\{ x\in \Omega;\, V(x)=m_0\}.
  \end{equation}
\end{itemize}
For $x_0\in K$, we have a formal limit problem:
    \begin{equation} \label{eq:1.8}
    -\Delta U + m_0U =  (I_\alpha * F(U)) f(U) \ \text{in $\R^N$}, \ U > 0, \ U \in H^1(\R^N),
    \end{equation}
whose solutions are characterized as critical points of 
    $$  L_{m_0}(u)=\half\int_{\R^N} |\nabla u|^2 + m_0u^2\, dx -\half \int_{\R^N} (I_\alpha*F(u))F(u)\, dx.
    $$
The nonlinear Choquard equation \eqref{eq:1.8} has been introduced and studied by Ma and Zhao in \cite{MZ} and successively in
\cite{MS,MS1} (see also the recent survey \cite{MS3}).
In particular, in \cite{MS1}, Moroz and Van Schaftingen proved that under (f1)--(f5) the limiting equation \eqref{eq:1.8} has
a positive ground state solution, which is radially symmetric about some point
(see also \cite{CCS,CS}).

However, differently from $(\ref{eq:1.4})$,
we do not have, in general, uniqueness and non-degeneracy for the ground state solutions of the limit equation \eqref{eq:1.8}
and it seems difficult
to perform the Lyapunov-Schmidt reduction method for deriving semiclassical states of  
\eqref{eq:1.5}.

In \cite{MS2},
Moroz and  Van Schaftingen proved existence of a single-peak solution for $(\ref{eq:1.1})$ concentrating at a local minima of $V$, 
when
$f(s)= |s|^{p-2} s$, $p \in [2, {N+\alpha\over N-2})$ 
using a novel nonlocal penalization technique,
        and recently in  \cite{YZZ} Yang, Zhang and Zhang showed the existence under (f4),
$\lim_{t\to\infty} {f(t)\over t^{\alpha+2\over N-2}} =0$ and
    \begin{equation}\label{eq:1.9}
    \lim_{t\to 0^+} {f(t)\over t} =0.
    \end{equation}
They also proved the existence of multi-peak solutions, whose each peak concentrates at different local minimum point
of $V(x)$ as $\e\to 0$.
Thus, for a singular perturbation problem, the situation is more delicate
and stronger conditions on the behavior of $f(s)$ as $s\sim 0$
are required in \cite{MS2,YZZ}.

We improve these results.  Our main result is

\medskip

\begin{theorem}\label{S:1.1}
  Suppose $N\geq 3$, (f1)--(f5), (V1)--(V2) and let $K$ be given in \eqref{eq:1.7}.
  Then \eqref{eq:1.1} has at least one
  positive solution concentrating in $K$ as $\e\to 0$.
\end{theorem}

As a special case, we have

\begin{corollary}\label{S:1.2}
  Suppose $N\geq 3$ and $f(s)=|s|^{p-2}s$ with $p\in ({N+\alpha\over N}, {N+\alpha\over N-2})$.  Moreover assume (V1)--(V2).
  Then \eqref{eq:1.1} has has at least one
  positive solution concentrating in $K$ as $\e\to 0$.
\end{corollary}

\begin{remark}
\begin{itemize}
\item[(i)]
  Our Theorem \ref{S:1.1} and Corollary \ref{S:1.2} can be applied to the locally sublinear case, e.g.,
  $f(s)=|s|^{p-2}s$, $p\in ({N+\alpha\over N},2)$.  In particular, it gives a positive answer to an open
  problem arisen in Moroz and Van Schaftingen \cite[Section 6.1]{MS2}.
\item[(ii)]
  For $f(s)=|s|^{p-2}s$, the range $({N+\alpha\over N}, {N+\alpha\over N-2})$ of $p$ is optimal for the existence of nontrivial solutions.
  In fact, in \cite[Theorem 2]{MS}. Moroz and Van Schaftingen showed that there are not nontrivial solutions outside.
\end{itemize}
\end{remark}

\medskip

To show our Theorem \ref{S:1.1} and Corollary \ref{S:1.2}, we develop a new variational argument  which relies and extends some works by Byeon and Jeanjean \cite{BJ} and Byeon and Tanaka \cite{BT1,BT2}.  
We remark that in \cite{BJ,BT1,BT2} they deal with singular perturbation problems for 
nonlinear Schr\"odinger equations, while our argument also works for the non-local problem \eqref{eq:1.1}.

Let $E(m_0)$ be the least energy level for the limit problem:
    $$  E(m_0)=\inf\{ L_{m_0}(\omega);\, \omega\in H^1(\R^N)\setminus\{ 0\},\, L_{m_0}'(\omega)=0\}
    $$
and we set for $s_0\in (0,\half)$
    \bea
    S_{m_0} &=& \{ \omega\in H^1(\R^N)\setminus\{ 0\};\, L_{m_0}'(\omega)=0,\, L_{m_0}(\omega)=E(m_0),\\ 
    &&\quad  \omega(0)=\max_{x\in\R^N}\omega(x)\}, \\
    \whS &=&\{ \omega(x/s);\, \omega\in S_{m_0},\, s\in [1-s_0,1+s_0]\}.
    \eea
We note that $S_{m_0}$ is compact.   We also set for $\nu>0$ small
    $$  N_\e(\nu)=\{ u(x)=U(x)+\varphi(x);\, U\in \whS,\, \|\varphi\|_{H^1} <\nu,\, \e\beta(u)\in K_{d_0}\},
    $$
where $\beta(u)$ is a center of mass of $u$, defined in Section \ref{section:3} and $K_{d_0}$ is a neighborhood of $K$
and we try to find critical points of $J_\e(u)$ in $N_\e(\nu)$.  We note that the parameter $s$ is introduced
corresponding to the scaling property of $L_{m_0}(u)$ and our definition of $\whS$, $N_\e(\nu)$ is slightly different
from those in \cite{BJ,BT1,BT2}. 

Such a localized variational method is used together with a penalization method in \cite{BJ,CJT,CJT2}
or a tail minimizing operator in \cite{BT1,BT2} to ensure that concentration of solutions occurs in 
the desired region.  In \cite{BJ,CJT,CJT2}, it is  introduced a penalization term and first a critical point of the penalized functional is found.

In \cite{BT1,BT2}, for a given $u\in N_\e(\nu)$ it is considered
a minimizing problem of type:
    \begin{eqnarray}
       \hbox{minimize}\, \{ J_\e(v);\, &&v\in H^1(\R^N)
        ,\, v(x)=u(x), \hbox{for}\, |x-\beta(u)|\geq R, \nonumber \\
                            &&\| v\|_{H^1(|x-\beta(u)|\geq R)} <\nu \} \ (R\gg 1) \label{eq:1.10}  
    \end{eqnarray} 
and it is defined the tail minimizing operator as a unique minimizer of \eqref{eq:1.10}, which helps to control
the behavior of solutions outside small balls.
Indeed, to justify such approaches, in \cite{BJ,CJT,CJT2} it is necessary to show that
solutions of penalized problems are small in the penalized region and this fact relies on the decay property of linear equations outside small balls.  To this aim it is crucial the exponential decay of solutions, which are valid if $f(s)=|s|^{p-2}s$,
$p\in [2,{N+\alpha\over N-2})$ or \eqref{eq:1.9} holds.
Of course, for the argument in \cite{BT1,BT2}, solvability and uniqueness of solutions of the minimizing
problem \eqref{eq:1.10} is required.

Moreover we note that $J_\e(u)$ defined in \eqref{eq:1.6} is well-defined and of class $C^1$ on $H^1(\R^N)$ 
under the conditions (f1)--(f5), especially for
$f(s)=|s|^{p-2}s$, $p\in ({N+\alpha\over N},{N+\alpha\over N-2})$.  However for a given $u\in H^1(\R^N)$, the minimizing
problem \eqref{eq:1.10} (or corresponding linear problem) outside small balls is not solvable 
due to non-local feature of the problem under (f1)--(f5) and in \cite{MS2,YZZ}, the authors assume
$f(s)=|s|^{p-2}s$, $p\in [2,{N+\alpha\over N-2})$ or \eqref{eq:1.9}.  Actually in \cite{MS2}, cases $p\in (2,{N+\alpha\over N-2})$ and
$p=2$ are discussed separately.  See \cite{MS2}, especially Section 6.1.

In this paper we take a new approach, which is inspired by \cite{BT1,BT2}, and we find an operator $\tau_\e$,
which has similar properties to the tail minimizing operator in \cite{BT1,BT2}, via a deformation argument.  
Our operator $\tau_\e:\, N_\e(\nu)\to N_\e(5\nu)$ has the following properties:
    \bea
    &&\tau_\e(u)(x)=u(x) \quad \hbox{for}\ |x-\beta(u)|\leq {1\over\sqrt\e},\\
    &&J_\e(\tau_\e(u)) \leq J_\e(u),\\
    &&\|\tau_\e(u)\|_{H^1(|x-\beta(u)|\geq {2\over\sqrt\e})} \leq \widetilde\rho_\e,
    \eea
where $\widetilde\rho_\e$ is independent of $u\in N_\e(\nu)$ and $\widetilde\rho_\e\to 0$ as $\e\to 0$.

Using the operator $\tau_\e$ together with the standard deformation argument, 
we can directly deal with $J_\e(u)$ in $N_\e(\nu)$.  Virtue of our approach is the following;
first we can find a critical point without introducing penalization, which simplifies the argument, second,
we argue without linear equations outside small balls nor the uniform exponential decay of solutions, so we are able to 
deal with locally sublinear case, for example, $f(s)=|s|^{p-2}s$, $p\in ({N+\alpha\over 2},2)$.
We believe that our new approach is applicable to a wide class of singular perturbation problems.

We also study a multiplicity of positive single-peak solutions of \eqref{eq:1.1} under (f1)--(f5), relating the number of
positive solutions to the topology of the critical set $K$ of $V(x)$ in the spirit of well-known results
Bahri, Coron \cite{BCo} and  Benci, Cerami  \cite{BC1,BC2,BCP} for semilinear elliptic problems with Dirichlet boundary condition.

Analysis of the topological changes between two level sets of $J_\e(u)$ in a small neighborhood $N_\e(\nu)$ of expected solutions
is essential in our argument.  This approach is taken for nonlinear Schr\"odinger equation and nonlinear magnetic Schr\"odinger
equations in \cite{CJT,CJT2}.  We extend these ideas to non-local setting.  We note that we do not require any monotonicity
assumption on $f$.  As a consequence we can not use advantage of Nehari manifolds and we use relative category and cup-length to
estimate the topological changes.

Our result is

\begin{theorem}\label{S:1.4}
  Suppose $N\geq 3$ and that (f1)--(f5) and (V1)--(V2) hold.  Then, for sufficiently small $\epsilon >0$, $(\ref{eq:1.1})$ has at least
  $\cuplength(K)+1$ positive solutions $v_{\e}^i$, $i=1, \dots, \cuplength(K)+1$ concentrating  as $\e \to 0 $ in $K$,
  where $\cuplength(K)$ denotes the cup-length defined with Alexander-Spanier cohomology with coefficients
  in the field $\F$.
\end{theorem}

\begin{remark}\label{S:1.5}
  If $K=\mathbb{S}^{N-1}$, the $N-1$ dimensional sphere in $\R^N$, then $\cuplength(K) + 1=  \cat(K) =2$.  If $K=\mathbb{T}^N$ is the $N$-dimensional torus, then $\cuplength(K) + 1 = \cat(K)= N + 1$.
  However in general $\cuplength(K) +1 \leq \cat(K)$.
\end{remark}

\begin{remark}\label{S:1.6}
  When we say that the solutions $v_{\e}^i$, $i=1, \dots, \cuplength(K)+1$ of Theorem \ref{S:1.4}
  ($v_\e$ of Theorem \ref{S:1.1} resp.) concentrate when $\e \to 0$ in K, we mean that there exists
  a maximum point $x_{\e}^i$ of $v_{\e}^i$ ($x_\e$ of $v_\e$ resp.)
  such that $\lim_{\e \to 0}dist(x_{\e}^i, K) =0$ ($\lim_{\e \to 0}dist(x_{\e}, K) =0$ resp.)
  and that for any such $x_{\e}^i$, $w_{\e}^i(x) = v_{\e}^i(\e(x+x_{\e}^i))$ ($w_{\e}(x) = v_{\e}(\e(x+x_{\e}))$ resp.)
  converges, up to a subsequence, uniformly to a least energy solution $U$ of
  \eqref{eq:1.8}.
\end{remark}


Since Theorem \ref{S:1.1} can be regarded as a special case of Theorem \ref{S:1.4},  we mainly deal with Theorem \ref{S:1.4}.
To study our non-local problem, we develop several new tools.  

This paper is organized as follows:
In Section \ref{section:2}, we develop useful estimates for a non-local term $\int_{\R^N} (I_\alpha*F(u))F(u)\, dx$,
which we use throughout this paper.  In Section \ref{section:3}, we study the limit problem \eqref{eq:1.8}.  As a typical
feature of \eqref{eq:1.8}, the problem has 3 different types of homogeneities and we employ new ideas as a new Pohozaev type
function $\Lambda(u)$.  This method can be useful to study other elliptic problem involving different types of scaling.
We also introduce a center of mass in $N_\e(\nu)$.   In Section \ref{section:4},we give a gradient estimate of $J_\e(u)$ in a neighborhood
of $\partial N_\e(\nu)$, which will play an essential role in the proof.  In Section \ref{section:5}, we develop a deformation
argument for a tail of $u\in N_\e(\nu)$, i.e., outside small balls, to construct a map $\tau_\e$.
In Sections \ref{section:6}--\ref{section:7}, using ideas from \cite{CJT,CJT2}, we analyze change of the topology of
level sets of $J_\e(u)$ in $N_\e(\nu)$ via relative cup-length and give a proof of Theorem \ref{S:1.4}.

For Theorem \ref{S:1.1}, we can give a simpler proof without using cup-length nor category, which we believe of interest itself.  
We give an outline of the simpler proof in Section \ref{section:8}, in which we use a localized mountain pass argument.


\setcounter{equation}{0}
\section{\label{section:2}  Preliminaries}
In what follows we use the notation:
\begin{eqnarray*}
    &&\norm u_{H^1} = \left(\intRN \abs{\nabla u}^2+u^2\, dx\right)^{1/2},\\
    &&\norm u_r = \left(\intRN \abs{u}^r\, dx\right)^{1/r} \quad \hbox{for} \ r\in [1,\infty),\\
    && B(p,R)=\{ x\in\R^N;\, |x-p| <R\},\\
    && B(p,R)^c=\R^N\setminus B(p,R) \quad \hbox{for}\ p\in \R^N\ \hbox{and}\ R>0.
\end{eqnarray*}
We study the multiplicity of solutions to $(\ref{eq:1.5})$ via a variational method. That is, we look for
critical points of the functional $J_\epsilon \in C^1(H^1(\R^N),\R)$ defined by
$$  J_\epsilon(u)=\half\norm{\nabla u}_2^2+\half\intRN V(\epsilon x)u^2\, dx
-\frac{1}{2} {\mathcal D}(u)
$$
where ${\mathcal D}(u) = \intRN (I_\alpha *F(u)) F(u)\, dx$.
The critical points of $J_\epsilon$ are clearly solutions of $(\ref{eq:1.5})$.


%
\begin{proposition}\label{S:2.1}
Let $p$, $r>1$ and $\alpha\in (0,N)$ with  $\frac 1 p + \frac 1 r = \frac{N+\alpha}N$. 
Then there exists a constant $C=C(N,\alpha,p,r)>0$ such that 
    \begin{equation} \label{eq:2.1}
    \left| \int_{\R^N} (I_\alpha*f) g\, dx \right| \leq C\|f\|_p \|g\|_{r}
    \end{equation}
for all $f\in L^p(\R^N)$ and $g\in L^r(\R^N)$.
\end{proposition}

Proposition \ref{S:2.1} can be derived from the following generalized Hardy-Littlewood-Sobolev inequality.
In the following proposition, $L_w^q(\mathbb{R}^N)$ is the \emph{weak} $L^q$ space.  
See \cite{Lieb1983,ll} for the definition.  
We denote by $\| \cdot\|_{q,w}$ the usual norm in $L_w^q(\mathbb{R}^N)$.

\begin{proposition}\label{S:2.2}
Let $p$, $q$, $r\in (1,\infty)$ satisfies $\frac  1 p +\frac 1 q +\frac 1 r =2$.
Then there exists a constant $N_{p,q,r}>0$ such that  for any 
$f \in L^p(\mathbb{R}^N)$, $g \in L^r(\mathbb{R}^N)$ and $h \in L_w^q(\mathbb{R}^N)$
    \begin{equation} \label{eq:2.2}
    \left| \int_{\R^N\times\R^N} f(x)h(x-y)g(y)\, dx\, dy
    \right| \leq N_{p,q,r} \|f\|_p \|g\|_r \|h\|_{q,w}.
    \end{equation}
\end{proposition}


In what follows, we denote  various constants, which are independent of $u\in H^1(\R^N)$, by $C$, $C'$, $C''$, $\cdots$.

By Proposition \ref{S:2.1}, we have for $u$, $v\in H^1(\R^N)$
    \bea
    |\calD(u)| &\leq& C\| F(u)\|_{2N\over N+\alpha}^2 
        \leq C\left( \|u\|_2^2 + \|u\|_{2N\over N-2}^{2N\over N-2}\right)^{N+\alpha\over N}, \\
    |\calD'(u)u| &\leq& C\| F(u)\|_{2N\over N+\alpha}\| f(u)u\|_{2N\over N+\alpha} 
        \leq C'\left( \|u\|_2^2 + \|u\|_{2N\over N-2}^{2N\over N-2}\right)^{N+\alpha\over N}, \\
    |\calD'(u)v| &\leq& C\| F(u)\|_{2N\over N+\alpha}\| f(u)v\|_{2N\over N+\alpha} \\
        &\leq& C'\left( \|u\|_2^2 + \|u\|_{2N\over N-2}^{2N\over N-2}\right)^{N+\alpha\over 2N}  
            \left( \|u\|_2^{2\alpha\over N+\alpha} + \|u\|_{2N\over N-2}^{{2N\over N-2}{2+\alpha \over N+\alpha}}\right)^{N+\alpha\over 2N}
            \| v\|_{H^1}.
    \eea
We also have
    \bea
    J_\e(u) &\geq& \half\|\nabla u\|_2^2 + {\underline V\over 2}\|u\|_2^2 -C\|F(u)\|_{2N\over N+\alpha}^2,\\
    J_\e'(u)u  &\geq& \|\nabla u\|_2^2 + {\underline V}\|u\|_2^2 
        -C\|F(u)\|_{2N\over N+\alpha} \|f(u)u\|_{2N\over N+\alpha}.
    \eea
In particular, $J_\e(u)$ has mountain pass geometry uniformly in $\e\in (0,1]$ and we have

\begin{corollary} \label{S:2.3}
There exist $\nu_0>0$ and $c_0>0$ such that
    $$  J_\e(u) \geq c_0\|u\|_{H^1}^2, \quad J_\e'(u)u \geq c_0\|u\|_{H^1}^2
    $$
for all $u\in H^1(\R^N)$ with $\|u\|_{H^1} \leq \nu_0$.
\end{corollary}


For latter uses, we give 

\begin{lemma}\label{S:2.4}
Let $p$, $r>1$ and $\alpha\in (0,N)$ with $\frac 1 p +\frac 1 r <\frac{N+\alpha}N$.  Then
there exists a constant $D_R>0$ depending on $R>0$ such that
    $$  D_R\to 0 \quad \text{as}\ R\to\infty
    $$
and
    $$  \left| \int_{\R^N} (I_\alpha*f)g \, dx\right| 
        \leq D_R \| f\|_p \|g\|_r
    $$
for all $f\in L^p(\R^N)$ and $g\in L^r(\R^N)$ with 
$\text{\rm dist}(\hbox{\rm supp}\, f,\hbox{\rm supp}\, g) \geq R$.
\end{lemma}

\medskip

\claim Proof.
We set
    $$  I_\alpha^R(x)=\begin{cases} 
            \frac 1 {|x|^{N-\alpha}} & \hbox{\rm if}\ |x| \geq R,\\
            0                       & \hbox{\rm if}\ |x| < R.\\
            \end{cases}
    $$
We note that $\|I_\alpha^R\|_{\frac N {N-\alpha},w}$ remains bounded as $R\to\infty$ but 
$\|I_\alpha^R\|_{\frac N {N-\alpha},w}\not\to 0$ as $R\to\infty$, while $\|I_\alpha^R\|_q\to 0$
as $R\to\infty$ for any $q\in (\frac N {N-\alpha}, \infty)$.

For $\frac 1 p +\frac 1 r <\frac{N+\alpha}N$, we find $q\in (\frac N {N-\alpha}, \infty)$ such
that $\frac 1 p +\frac 1 q +\frac 1 r =2$.  For $f\in L^p(\R^N)$ and $g\in L^r(\R^N)$ with
$\text{\rm dist}(\hbox{\rm supp}\, f,\hbox{\rm supp}\, g) \geq R$,
the Hausdorff-Young inequality implies
    $$
    \left|\int_{\R^N} (I_\alpha*f)g \, dx\right| 
    = \left|\int_{\R^N} (I_\alpha^R*f)g \, dx\right|  
    \leq  \| I_\alpha^R\|_q \|f\|_p \|g\|_r.
    $$
Setting $D_R=\| I_\alpha^R\|_q$, we have the conclusion.  \QED

As a corollary to Proposition \ref{S:2.1} and Lemma \ref{S:2.4}, we have

\begin{corollary}\label{S:2.5}
Suppose that $F(s):\,\R\to\R$ satisfy (f1)--(f3).  Then there exists a
constant $m_R>0$ such that 
    $$  m_R\to 0 \quad \hbox{as}\ R\to \infty
    $$
and for any $u\in H^1(\R^N)$ and $g\in L^{\frac{2N}{N+\alpha}}(\R^N)$
with $\text{\rm dist}(\hbox{\rm supp}\, u,\hbox{\rm supp}\, g) \geq R$
    $$  \left|\int_{\R^N}(I_\alpha*F(u)) g \,dx \right|
        \leq m_R \sigma(\|u\|_{H^1})^{\frac{N+\alpha}{2N}}
            \|g\|_{\frac{2N}{N+\alpha}},
    $$
where
    \begin{equation}\label{eq:2.3}
    \sigma(t)=t^2 + t^{2N\over N-2}.
    \end{equation}
\end{corollary}

\claim Proof.
We fix $q\in (2,\frac{2N}{N-2})$ and write $q={2N\over N-2}\theta$ with $\theta\in ({N-2\over N},1)$.  
By (f3), for any $\delta>0$ there exists $C_\delta>0$ such
that
    $$  |F(s)| \leq \delta \sigma(s)^{N+\alpha\over 2N} + C_\delta|s|^{q{N+\alpha\over 2N}}
        \quad \hbox{for all}\ s\in\R.
    $$
Setting
    \begin{eqnarray*}
    H_1(s)&=&\begin{cases}
        \delta \sigma(s)^{N+\alpha\over 2N}
            & \hbox{if}\ |F(s)| \geq \delta \sigma(s)^{N+\alpha\over 2N}, \\
        |F(s)|  & \hbox{if}\ |F(s)| < \delta \sigma(s)^{N+\alpha\over 2N}, 
        \end{cases}\\
    H_2(s)&=&|F(s)|- H_1(s),
    \end{eqnarray*}
we have
    \begin{eqnarray*}
    |F(s)| &=& H_1(s) +H_2(s), \\
    H_1(s) &\leq& \delta \sigma(s)^{N+\alpha\over 2N}, \\
    H_2(s) &\leq& C_\delta |s|^{q{N+\alpha\over 2N}}.
    \end{eqnarray*}
We note for $u\in H^1(\R^N)$
    \bea
    \| H_1(u)\|_{2N\over N+\alpha} &\leq& \delta \left( \|u\|_2^2 +\|u\|_{2N\over N-2}^{2N\over N-2}\right)^{N+\alpha\over 2N}
        \leq \delta C' \sigma(\|u\|_{H^1})^{N+\alpha\over 2N},\\
    \| H_2(u)\|_{{2N\over N+\alpha}{1\over\theta}} &\leq& C_\delta \|u\|_{q/\theta}^{q{N+\alpha\over 2N}} 
        = C_\delta \|u\|_{2N\over N-2}^{q{N+\alpha\over 2N}}  
        \leq C_\delta' \|u\|_{H^1}^{q{N+\alpha\over 2N}} \\
        &\leq& C_\delta'' \sigma(\|u\|_{H^1})^{N+\alpha\over 2N}.
    \eea
Thus, 
    \bea
    \left|\int_{\R^N} (I_\alpha*F(u))g\,dx\right|
    &\leq& \int_{\R^N} (I_\alpha*H_1(u))|g|\,dx + \int_{\R^N} (I_\alpha*H_2(u))|g|\,dx \\
    &\leq& C\|H_1(u)\|_{2N\over N+\alpha}\|g\|_{2N\over N+\alpha}
            + D_R\|H_2(u)\|_{{2N\over N+\alpha}{1\over\theta}} \|g\|_{2N\over N+\alpha},\\
    &\leq& (CC'\delta + D_R C_\delta'')
        \sigma(\|u\|_{H^1})^{\frac{N+\alpha}{2N}}
            \|g\|_{\frac{2N}{N+\alpha}}.
    \eea
We set
    \begin{eqnarray*}   
    m_R&\equiv& \sup_{u\in H^1(\R^N)\setminus\{0\}} 
    {\left|\int_{\R^N}(I_\alpha*F(u)) g \,dx \right| \over
        \sigma(\|u\|_{H^1})^{\frac{N+\alpha}{2N}}
            \|g\|_{\frac{2N}{N+\alpha}}} \\
    &\leq& CC'\delta + D_R C_\delta''.
    \end{eqnarray*}
It follows from Lemma \ref{S:2.4} that
    $$  \limsup_{R\to\infty} m_R \leq CC'\delta.
    $$
Since $\delta>0$ is arbitrary, we have $m_R\to 0$ as $R\to\infty$.
\QED

\medskip

For $R>0$ we choose functions $\zeta_R(s)$, $\widetilde\zeta_R(s)\in C^\infty(\R,\R)$ such that
    \begin{equation}\label{eq:2.4}
    \zeta_R(s)=\begin{cases} 1 &\hbox{for $s\leq R$},\\ 0 &\hbox{for $s\geq R+2$}, \end{cases} \quad
    \widetilde \zeta_R(s)=\begin{cases} 0 &\hbox{for $s\leq R-2$},\\ 1 &\hbox{for $s\geq R$}, \end{cases} 
    \end{equation}
and $\zeta_R(s)$, $\widetilde\zeta_R(s)\in [0,1]$, $-\zeta_R'(s)$, $\widetilde\zeta_R'(s)\in [0,1]$ for all $s\in\R$.
We also set
    $$  \chi_R(s)=\begin{cases} 1 &\hbox{for $s\leq R$}, \\ 0 & \hbox{for $s>R$} \end{cases}, \quad 
        \widetilde\chi_R(s)=1-\chi_R(s).
    $$

\begin{corollary}\label{S:2.6}
For a fixed $M>0$, suppose that $u$, $v\in H^1(\R^N)$ satisfy
    $$  \|u\|_{H^1}, \ \|v\|_{H^1} \leq M.
    $$
Then there exist $C_1$, $C_2>0$ such that for any $p\in \R^N$ and $R$, $L>4$ we have
\begin{itemize}
\item[(i)] $|\calD(u)-\calD(\zeta_R(|x-p|) u)-\calD(\widetilde \zeta_{R+L}(|x-p|)u)| 
\leq C_1\|u\|_{H^1(|x-p|\in [R,R+L])}^{N+\alpha\over N}
+C_2 m_L$.
\item[(ii)] $|\calD'(u)v -\calD'(\zeta_R(|x-p|) u)(\zeta_R(|x-p|)v)  \hfil\break
-\calD'(\widetilde \zeta_{R+L}(|x-p|)u)(\widetilde \zeta_{R+L}(|x-p|)v)| \hfil\break
\leq C_1 \max\{ \|u\|_{H^1(|x-p|\in [R,R+L])}, \|v\|_{H^1(|x-p|\in [R,R+L])} \}^{N+\alpha\over N} +C_2m_L$.
\item[(iii)] $| (\calD'(u) - \calD'(\widetilde \zeta_{R+L}(|x-p|)u))(\widetilde\zeta_{R+L}(|x-p|)v) | \\
    \leq C_1 \max\{ \|u\|_{H^1(|x-p|\in [R,R+L])}, \|v\|_{H^1(|x-p|\in [R,R+L])} \}^{N+\alpha\over N} +C_2m_L$.
\end{itemize}
\end{corollary}

\claim Proof.
We may assume $p=0$ and we prove just (ii).  (i) and (iii) can be shown in a similar way.

We note that 
    $$  \calD'(u)v =\int_{\R^N} (I_\alpha*F(u)) f(u)v\, dx.
    $$
We set 
    \bea
    &&F_A=F(\chi_R(|x|) u), \quad F_B=F(\widetilde \chi_{R+L}(|x|)u), \\
    &&F_C=F((1-\chi_R(|x|)-\widetilde \chi_{R+L}(|x|))u),\\
    &&g_A=f(\chi_R(|x|) u)\chi_R(|x|)v, \quad g_B=f(\widetilde \chi_{R+L}(|x|) u)\widetilde\chi_{R+L}v(|x|), \\
    &&g_C= f((1-\chi_R(|x|)-\widetilde \chi_{R+L}(|x|))u)((1-\chi_R(|x|)-\widetilde \chi_{R+L}(|x|))v).  
    \eea
We also set
$\delta=\max\{ \|u\|_{H^1(|x|\in [R,R+L])}, \|v\|_{H^1(|x|\in [R,R+L])} \}$.

First we show for some constants $C_1$, $C_2>0$ independent of $u$, $v$
    \begin{eqnarray}
    &&\left| \int_{\R^N}(I_\alpha*F(u))f(u)v\, dx - \int_{\R^N}(I_\alpha*F_A)g_A\, dx  - \int_{\R^N}(I_\alpha*F_B))g_B\, dx\right| \nonumber\\
    &\leq& C_1\delta^{N+\alpha\over N} +C_2m_L. \label{eq:2.5}
    \end{eqnarray}
In fact,
    \begin{eqnarray*}
    &&\int_{\R^N}(I_\alpha*F(u))v\, dx \\
    &=& \int_{\R^N} (I_\alpha*(F_A+F_B+F_C)) (g_A+g_B+g_C)\,dx \\
    &=& \int_{\R^N} (I_\alpha*F_A) g_A\, dx + \int_{\R^N} (I_\alpha*F_B) g_B \, dx\\
    && +\int_{\R^N} (I_\alpha*F_A) g_B\, dx + \int_{\R^N} (I_\alpha*F_B) g_A \, dx\\
    && +\int_{\R^N} (I_\alpha*F_C) (g_A+g_B)\, dx + \int_{\R^N} (I_\alpha*(F_A+F_B)) g_C\, dx \\
    && +\int_{\R^N} (I_\alpha*F_C) g_C \, dx
    \end{eqnarray*}
By Corollary \ref{S:2.5},
    $$  \left|\int_{\R^N} (I_\alpha*F_A) g_B\, dx + \int_{\R^N} (I_\alpha*F_B) g_A \, dx\right|
        \leq Cm_L.
    $$
We can easily see that for some constant $C$ depending only on $M$
    \begin{eqnarray*}
    &&\|F_A\|_{2N\over N+\alpha}, \ \|F_B\|_{2N\over N+\alpha}, \ \|g_A\|_{2N\over N+\alpha}, \ \|g_B\|_{2N\over N+\alpha}
        \leq C, \\
    &&\|F_C\|_{2N\over N+\alpha}, \ \|g_C\|_{2N\over N+\alpha} \leq C\delta^{N+\alpha\over N}.
    \end{eqnarray*}
Thus we get \eqref{eq:2.5}.  

Since
    \begin{eqnarray*}   
        &&\|F_A-F(\zeta_Ru)\|_{2N\over N+\alpha},\ \|F_B-F(\widetilde\zeta_{R+L}u)\|_{2N\over N+\alpha}, \\
        &&\|g_A-f(\zeta_Ru)(\zeta_Rv)\|_{2N\over N+\alpha},\ \|g_B-f(\widetilde\zeta_{R+L}u)(\widetilde\zeta_{R+L}v)\|_{2N\over N+\alpha}
        \leq C\delta^{N+\alpha\over \alpha},
    \end{eqnarray*}
we have the conclusion (ii).  \QED

\medskip
    
We will use Corollary \ref{S:2.6} with $M=\sup_{\omega\in S_{m_0},\, s\in [\half,{3\over 2}]} 
\| \omega(x/s)\|_{H^1}+1$ repeatedly.


\setcounter{equation}{0}
\section{\label{section:3} Limit problem}
\subsection{\label{subsection:3.1} Limit problem}

Let $m_0$, $K$ be introduced in (V2).  For $d>0$, let $K_d=\{ x\in \R^N;\, \inf_{y\in K} |x-y| <d\}$ be
a $d$-neighborhood of $K$.
    
In what follows we denote 
$K=\{ x\in\Omega ; \, V(x)=m_0\}$ where $m_0$ is introduced in $\hbox{(V1)}$ and  $K_d$ a $d$-neighborhood of $K$.

We can choose $d_0>0$ small such that $K_{2d_0}\subset \Omega$ and
\begin{equation}\label{eq:3.1}
V(x) \geq m_0 + \rho_0 \quad \hbox{for all}\  x \in K_{2d_0}\setminus K_{\half d_0}
\end{equation}
for some $\rho_0>0$ and the conclusion of Lemma \ref{S:7.5} in the following Section \ref{section:7} holds for
$d=d_0$.

\smallskip
For $a>0$ we define a functional $L_a \in C^1(H^1(\R^N),\R)$ by
    $$  L_a(u)=\half\norm{\nabla u}_2^2+{a\over 2}\norm u_2^2 - \frac{1}{2} {\mathcal D}(u).
    $$
Especially $L_{m_0}(u)$ is associated to the limit problem
\begin{equation}\label{eq:3.2}
- \Delta u + m_0 u = (I_\alpha * F(u)) f(u), \quad u \in H^1(\R^N).
\end{equation}

\smallskip

We denote by $E(m_0)$ the least energy level for \eqref{eq:3.2}.
That is,
$$  
E(m_0)=\inf\{ L_{m_0}(u);\, u\not=0,\ L_{m_0}'(u)=0\}.
$$
In \cite{MS1} it is proved that there exists a least energy solution of \eqref{eq:3.2}
if (f1)--(f4) are satisfied.
They proved that if $f$ is odd and has a constant sign on $(0, +\infty)$, then every ground state solution of \eqref{eq:3.2} 
has constant sign and is radially symmetric with respect to some point in $\R^N$.

Also it is showed that each solution of \eqref{eq:3.2} satisfies the Pohozaev's identity
\begin{equation}\label{eq:3.3}
\frac{N-2}{2} \|\nabla u\|_2^2 + {N \over 2}m_0 \|u\|_2^2 = \frac{N + \alpha}{2} \, {\mathcal D}(u).
\end{equation}

We set
$$S_{m_0} = \{  \omega \in H^1(\R^N)\setminus\{ 0\} ; \,
L_{m_0}'(\omega)=0,\ L_{m_0}(\omega) =E(m_0), \ \omega(0)=\max_{x\in\R^N} \omega(x)\}.
$$

Arguing as in \cite{MS1} we can prove that $S_{m_0}$  is compact in $H^1(\R^N)$ and that its elements have a uniform
decay.  Especially,
    \begin{equation}\label{eq:3.4}  
    \|\omega\|_{H^1(B(0,R)^c)} \to 0 \quad \hbox{as}\ R\to\infty \ \hbox{uniformly in}\ \omega\in S_{m_0}.
    \end{equation}
Moreover, they have a uniform exponential decay if $\limsup_{s\to 0}{f(x)\over s}<\infty$ and they have a uniform
polynomial decay if $F(s)=|s|^p$ with $p\in ({N+\alpha\over N},2)$.  See Moroz-Van Schaftingen \cite{MS}.

\subsection{\label{subsection:3.2} A  Pohozaev type function}
To study the scaling property of the limit equation \eqref{eq:3.2}, we define a functional
$$
J(\lambda, u) = L_{m_0}(u({\cdot}/{\lambda}))= \frac{\lambda^{N-2}}{2} \norm{\nabla u}_2^2+ m_0 \frac{\lambda^{N}}{2}  \norm u_2^2
- \frac{\lambda^{N + \alpha}}{2} {\mathcal D}(u).
$$
For any $u\in H^1(\R^N)\setminus \{ 0\}$, we have
$$
\frac{\partial}{\partial \lambda}  J(\lambda, u) =\frac{(N-2) \lambda^{N-3}}{2}  \norm{\nabla u}_2^2+ m_0 \frac{N \lambda^{N-1}}{2}  \norm u_2^2
- \frac{(N+ \alpha) \lambda^{N + \alpha - 1}}{2} {\mathcal D}(u).
$$
Since $\lim_{\lambda\to 0^+}{1\over \lambda^{N-2}} J(\lambda,u)>0$, $\lim_{\lambda\to\infty} J(\lambda,u)=-\infty$ and
moreover for any critical point $\lambda\in (0,\infty)$ of $\lambda\mapsto J(\lambda,u)$,
    $$  {\partial^2\over\partial\lambda^2}J(\lambda,u)
        = - \left\{ {(\alpha+2)(N-2)\over 2} \lambda^{N-4} \|\nabla u\|_2^2
            + {\alpha Nm_0\over 2}\lambda^{N-2}\|u\|_2^2 \right\} <0,
    $$
we observe that for any $u\in H^1(\R^N)\setminus\{ 0\}$, $(0,\infty)\to\R;\, \lambda\mapsto J(\lambda,u)$
has a unique critical point $\Lambda(u)$, which is a non-degenerate local maximum.
We note that $\Lambda(u):\, H^1(\R^N)\setminus\{ 0\}\to\R$ is continuous and by the Pohozaev
identity,
    \begin{equation}\label{eq:3.5}  \Lambda(u)=1 \quad \hbox{for all}\ u\in S_{m_0}.
    \end{equation}
We also have
    \begin{equation}\label{eq:3.6}
    \Lambda(u)=s \quad \hbox{for}\ u=\omega({x-p\over s})\ \hbox{with}\ \omega\in S_{m_0},\ s>0,\ p\in\R^N.
    \end{equation}

\begin{lemma}\label{S:3.1}
There exist $C_0>0$ and $s_0\in (0,\half)$ such that for any $v\in H^1(\R^N)\setminus\{ 0\}$ with $\Lambda(v)=1$
and $s\in [1-2s_0,1+2s_0]$
    $$  L_{m_0}(v({x\over s})) \geq E(m_0)(1-C_0(s-1)^2).
    $$
\end{lemma}

\claim Proof.
As in \cite{MS1} (see also \cite{JT1}), we can prove the following characterization:
    \be \label{eq:3.7}
    E(m_0)=\inf\{ L_{m_0}(u);\, u\not=0,\, \Lambda(u)=1\}.
    \ee
For $v\in H^1(\R^N)\setminus\{ 0\}$ with $\Lambda(v)=1$ and $s>0$,
    \begin{eqnarray*}
    &&L_{m_0}(v({x\over s})) = {s^{N-2}\over 2}\|\nabla v\|_2^2 +{s^N\over 2}m_0\| v\|_2^2
            -{s^{N+\alpha}\over 2}\calD(v) \\
    &&\quad = \half (s^{N-2}-{N-2\over N+\alpha}s^{N+\alpha}) \|\nabla v\|_2^2
        + \half(s^N-{N\over N+\alpha}s^{N+\alpha}) m_0\|v\|_2^2 \\
    &&\quad\equiv \half g(s) \|\nabla v\|_2^2 +\half h(s) m_0\| v\|_2^2.
    \end{eqnarray*}
It is easy to see that $g(s)$ and $h(s)$ take their unique maxima at $s=1$ and $g''(1)<0$, $h''(1)<0$.
Thus there exists $C_0>0$ and $s_0\in (0,\half)$ such that
    \begin{eqnarray*}
    g(s) &\geq& g(1)(1-C_0(s-1)^2), \\
    h(s) &\geq& h(1)(1-C_0(s-1)^2)
    \end{eqnarray*}
for all $s\in [1-2s_0,1+2s_0]$.  Therefore we have
    \begin{eqnarray*}
    L_{m_0}(v({x\over s})) &\geq& \left(\half g(1)\|\nabla v\|_2^2+\half h(1)m_0\| v\|_2^2\right) (1-C_0(s-1)^2) \\
    &=& L_{m_0}(v) (1-C_0(s-1)^2).
    \end{eqnarray*}
By \eqref{eq:3.7}, we have the conclusion of Lemma \ref{S:3.1}. \QED

\medskip

\begin{corollary}\label{S:3.2}
Assume that $u\in H^1(\R^N)\setminus\{ 0\}$ satisfies $\Lambda(u)\in [1-2s_0,1+2s_0]$.  Then
    $$  L_{m_0}(u) \geq E(m_0) (1-C_0(\Lambda(u)-1)^2).
    $$
\end{corollary}

\claim Proof.
Set $v(x)=u(\Lambda(u)x)$.  Then we have $\Lambda(v)=1$ and $v(x)=v({x\over \Lambda(u)})$.  Thus Corollary \ref{S:3.2}
follows from Lemma \ref{S:3.1}.  \QED

\medskip

\begin{corollary}\label{S:3.3}
Choosing $s_0$ smaller, there exists  $\delta_0>0$ such that for all $\omega\in S_{m_0}$
and $s\in [1-2s_0, 1+2s_0]$
    \be\label{eq:3.8}
    L_{m_0+\rho_0}(\omega({x\over s})) \geq E(m_0) + 3\delta_0.
    \ee
\end{corollary}

\claim Proof.
Let $s_0\in (0,\half)$ be given in Lemma \ref{S:3.1}.  
For $\omega\in S_{m_0}$ and $s\in [1-2s_0, 1+2s_0]$ we have
    \bea    
    &&L_{m_0+\rho_0}(\omega({x\over s})) = L_{m_0}(\omega({x\over s})) + \half \rho_0 \| \omega({x\over s})\|_2^2 \\
    &&\quad \geq E(m_0)(1-C_0(s-1)^2) +\half \rho_0 \min_{\omega\in S_{m_0},\, s\in [1-2s_0, 1+2s_0]} \| \omega({x\over s})\|_2^2.
    \eea
Setting $\delta_0={1\over 6} \rho_0 \min_{\omega\in S_{m_0},\, s\in [1-2s_0, 1+2s_0]} \| \omega({x\over s})\|_2^2>0$,
and choosing a smaller $s_0$,  we have \eqref{eq:3.8}.  \QED

\medskip

In what follows, we fix $s_0\in (0,\half)$ for which the conclusion of Lemma \ref{S:3.1} and \eqref{eq:3.8} holds.


\subsection{\label{subsection:3.3} A center of mass}
Let $s_0\in (0,\half)$ be the number given at the end of the previous section.  We set
    $$  \whS = \{  \omega({x-p\over s});\, \omega\in S_{m_0},\ p\in\R^N,\ s\in [1-s_0,1+s_0]\}.
    $$
We also set
    \bea
    \wrho(u) &=& \inf_{ U\in \whS} \, \|u-U\|_{H^1} \\
    &=& \inf\{ \| u- \omega({x-p\over s})\|_{H^1};\,  \omega\in S_{m_0},\ p\in\R^N,\ s\in [1-s_0,1+s_0]\}.
    \eea
For $\nu>0$ small we denote $\nu$-neighborhood of $\whS$ by $\whS(\nu)$:
    \bea
    &&\whS(\nu) = \{ u\in H^1(\R^N);\, \wrho(u) <\nu\} \\
    &&= \{ \omega({x-p\over s})+\varphi(x);\, \omega\in S_{m_0},\, p\in\R^N,\ s\in [1-s_0,1+s_0], \, \|\varphi\|_{H^1} <\nu\}.
    \eea
By the compactness of $S_{m_0}$ and \eqref{eq:3.5}, we have

\begin{lemma}\label{S:3.4}
There exists $\nu_1>0$ such that for all $u\in \whS(\nu_1)$
    $$  \Lambda(u)\in [1-2s_0,1+2s_0].
    $$
\end{lemma}

\noindent
Following \cite{BT1, BT2, CJT} for $\nu_2>0$ small we introduce a center of mass in $\nu_2$-neighborhood $\whS(\nu_2)$ of $\whS$.

\begin{lemma}\label{S:3.5}
    There exist $\nu_2>0$, $R_0>0$ and a map $\beta:\, \widehat S(\nu_2) \to \R^N$
    such that
    $$  \abs{\beta(u)-p} \leq R_0
    $$
    for all $u(x)=\omega (\frac{x-p}{s})+\varphi(x)\in \widehat S (\nu_2)$ with $p\in\R^N$, $\omega \in S_{m_0}$, $s \in [1-s_0, 1+ s_0]$,
    $\norm\varphi_{H^1}< \nu_2$.  Moreover, $\beta(u)$ has the following properties:
    \begin{itemize}
        \item[(i)] $\beta(u)$ is shift equivariant, that is,
        $$  \beta(u(x-y))=\beta(u(x))+y   \quad \hbox{for all}\ u \in\widehat S (\nu_2)\
        \hbox{and}\ y\in\R^N.
        $$
        \item[(ii)] $\beta(u)$ is locally Lipschitz continuous, that is, there exist
        constants $c_1$, $c_2>0$ such that
        $$  \abs{\beta(u)-\beta(v)} \leq c_1\norm{u-v}_{H^1}
        \quad \hbox{for all}\  u, v\in\widehat S (\nu_2) \  \hbox{with}\
        \norm{u-v}_{H^1}\leq c_2.
        $$
        \item[(iii)] $\beta( \omega(\frac{x-p}{s})) =p$ for all $p \in \R^N$, $\omega \in S_{m_0}$, $s \in [1-s_0, 1+ s_0]$.
    \end{itemize}
\end{lemma}
The proof is given in \cite{BT1, BT2} in a slightly different situation.
We give here a simple proof.

\medskip

\claim Proof.
Set $\widetilde S_{m_0}=\{U(x)=\omega({x\over s});\, \omega\in S_{m_0},\, s\in [1-s_0,1+s_0]\}$.  Then we have
$\whS=\{ U(x-p);\, U\in \widetilde S_{m_0},\, p\in \R^N\}$.

Taking into account that
$\widetilde S_{m_0}$ is compact and  the uniform decay 
(\ref{eq:3.4}), we can define  $r_*=\min_{U \in\widetilde S_{m_0}}\norm U_{H^1}>0$ and choose $R_*>1$ such that for $U \in\widetilde S_{m_0}$
$$  \norm U_{H^1(\abs x\leq R_*)} > {3\over 4}r_*  \quad \mbox{and} \quad
\norm U_{H^1(\abs x\geq R_*)} < {1\over 8}r_*.
$$
For $u\in H^1(\R^N)$ and $p\in\R^N$, we define
$$  d(p,u)=\psi\left(\inf_{U \in \widetilde S_{m_0}} \norm{u-U(x-p)}_{H^1(\abs{x-p}\leq R_*)}\right),
$$
where $\psi \in C_0^\infty(\R,\R)$ is such that
\begin{eqnarray*}
    &&\psi(r)=\begin{cases}    1   &r\in [0,{1\over 4}r_*],\\
        0   &r\in [\half r_*,\infty),\end{cases}\\
    &&\psi(r)\in [0,1]\quad \hbox{for all}\ r\in [0,\infty).
\end{eqnarray*}
Now let
$$  \beta(u)={\displaystyle \intRN q\, d(q,u)\, dq
    \over\displaystyle \intRN d(q,u)\, dq}
\quad \hbox{for}\quad u \in\widehat S({1\over 8}r_*).
$$
We shall show that  $\beta$ has the desired property.
\vspace{2mm}

Let $u\in \widehat S({1\over 8}r_*)$ and write 
$u(x)= U(x-p)+\varphi(x)$ with $U\in \widetilde S_{m_0}$, $p\in\R^N$, $\norm\varphi_{H^1}\leq {1\over 8}r_*$.
 
Then for $\abs{q-p}\geq 2R_*$ and $\tilde U\in\widetilde S_{m_0}$ we have 
\begin{eqnarray*}
    \norm{u-\tilde U(x-q)}_{H^1(\abs{x-q}\leq R_*)}
    &\geq& \norm{\tilde U(x-q)}_{H^1(\abs{x-q}\leq R_*)}\\
    &&    - \norm{U ({x-p})}_{H^1(\abs{x-p}\geq R_*)}-{1\over 8}r_*\\
    &>& {3\over 4}r_*-{1\over 8}r_*-{1\over 8}r_*=\half r_*.
\end{eqnarray*}
Thus $d(q,u)=0$ for $\abs{q-p}\geq 2R_*$.  We can also see that, for small $r>0$
$$  d(q,u)=1 \quad \hbox{for}\ \abs{q-p}<r.
$$
Thus $B(p,r)\subset \supp d(\cdot,u)\subset B(p,2R_*)$.
Therefore $\beta(u)$ is well-defined and we have
$$  \beta(u)\in B(p,2R_*) \quad \hbox{for}\quad  u\in \widehat S({1\over 8}r_*).
$$
Shift equivariance and locally Lipschitz continuity of $\beta$ can be
checked easily.  Setting $\nu_2={1\over 8}r_*$ and $R_0=2R_*$, 
we have the desired result.  \QED
\vspace{3mm}


\bigskip

Using this lemma we have

\begin{lemma}\label{S:3.6}
For $\delta_0>0$ given in Corollary \ref{S:3.3}, we have for $\epsilon>0$ small
    $$  J_\epsilon (\omega( \frac{x - p/ \epsilon}{s})) 
        \geq E(m_0) + 2\delta_0.
    $$
for all $\omega \in S_{m_0}$, $s \in [1- s_0, 1+ s_0]$, $p \in K_{{3\over 2}d_0}\setminus K_{\half d_0}$.
\end{lemma}

\claim Proof.
For $\omega \in S_{m_0}$, $s \in [1- s_0, 1+ s_0]$, we have
    \begin{eqnarray*}
    J_\epsilon(\omega({x - p/ \epsilon\over s})) &=& 
    L_{V(p)}(\omega(x/s)) + \frac{1}{2} \int_{\R^N} (V(\epsilon x + p) -V(p)) \omega(x/s)^2\, dx \\
    &\geq& L_{m_0+\rho_0}(\omega(x/s)) + o (1),
    \end{eqnarray*}
where $o(1)$ tends to zero as $\epsilon \to 0$.

By \eqref{eq:3.8}, we have the desired result for $\epsilon>0$ small.  \QED

\bigskip
It follows the corollary:
\smallskip

\begin{lemma}\label{S:3.7}
There exists $\nu_3>0$  such that for $\epsilon>0$ small
    $$  J_\epsilon(u) \geq E(m_0) +\delta_0
    $$
for all $u \in \whS(\nu_3)$ with $\e\beta(u)\in K_{{5\over 4}d_0}\setminus K_{{3\over 4}d_0}$.
\end{lemma}

\claim Proof.
For $u(x)=\omega({x-p/\e\over s})+\varphi(x)$ with $\|\varphi\|_{H^1}<\nu_2$, we have
    $$  |\e\beta(u)-p| \leq \e|\beta(u)-{p\over \e}| \leq \e R_0.
    $$
Thus $\e\beta(u)\to p$ as $\e\to 0$ and for $\e>0$ small, $\e\beta(u)\in K_{{5\over 4}d_0}\setminus K_{{3\over 4}d_0}$
implies $p\in K_{{3\over 2}d_0}\setminus K_{\half d_0}$.
Since
    $$  J_\e(u)-J_\e(\omega({x-p/\e\over s})) \to 0
    $$
uniformly as $\|\varphi\|_{H^1}\to 0$,  Lemma \ref{S:3.6} implies the conclusion of Lemma \ref{S:3.7} for small $\nu_3>0$.  \QED,


\setcounter{equation}{0}
\section{\label{section:4} Gradient estimates}
For $\nu>0$ small, we set
    $$  N_\e(\nu)=\{ u\in \whS(\nu);\, \e\beta(u)\in K_{d_0}\}.
    $$
We give the following $\e$-dependent concentration-compactness type result, which will give a useful gradient estimate
later in Corollary \ref{S:4.2}.

\begin{proposition}\label{S:4.1}
For sufficiently small $\nu_4>0$, $J_\e(u)$ has the following property:
Let $(\epsilon_j) \subset (0, 1)$ be a sequence such that $\epsilon_j \to 0$ as $j \to + \infty$
and let $(u_j) \subset H^1(\R^N)$ be such that $u_j\in N_{\epsilon_j}(\nu_4)$ and
\begin{eqnarray}
&&J_{\epsilon_j}(u_j) \to E(m_0),                       \label{eq:4.1}\\
&&J'_{\epsilon_j}(u_j) \to 0 \quad \hbox{as} \ j \to + \infty.   \label{eq:4.2}
\end{eqnarray}
Then there exist, up to a subsequence,   $(p_j) \subset \R^N$, $p_0 \in K$  and 
$\omega_0 \in S_{m_0}$ such that $p_j \to p_0$ and
$$
\| u_j\bigl(\cdot + \frac{p_j}{\epsilon_j}\bigr) -\omega_0 \|_{H^1} \to 0
\quad \hbox{as} \ j \to + \infty.
$$
\end{proposition}

\claim Proof.
Choose $\nu\in (0,\min\{ \nu_0,\nu_1,\nu_2, \nu_3\})$ such that
    $$  \nu < \half\min\{ \|\omega(x/s)\|_{H^1(B(0,R_0))};\, \omega\in S_{m_0},\, s\in [1-s_0,1+s_0]\}.
    $$
Suppose that $u_j\in N_{\e_j}(\nu)$ satisfies \eqref{eq:4.1}--\eqref{eq:4.2}.  Setting $p_j=\e_j\beta(u_j)\in K_{d_0}$, 
we have for some $p_0\in \overline{K_{d_0}}$ and $u_0\in H^1(\R^N)$
    \bea
    &&u_j(x+{p_j\over\e_j}) \wlimit u_0(x)\not=0 \quad \hbox{weakly in}\ H^1(\R^N),\\
    && p_j \to p_0\in \overline{K_{d_0}}.
    \eea
We will show that $u_j(x+{p_j\over\e_j})\to u_0(x)$ strongly in $H^1(\R^N)$ and $p_0\in K$.

\smallskip

\noindent{\sl Step 1:} For large $n\in\N$ and $j\in \N$, there exists $k_j^n\in \{ 1,2,\cdots,n\}$ such that
    $$  \| u_j\|_{H^1(A_{k_j^n}^n)}^2 \leq {4\nu^2\over n},
    $$
where
    $$  A_k^n =\{ x\in\R^N;\, nk \leq |x-{p_j\over \e_j}|\leq n(k+1)\}.
    $$

\smallskip

\noindent
In fact, for large $n$ we may assume $\| u_j\|_{H^1(B({p_j\over \e_j},n)^c)} \leq 2\nu$.  Since
    $$  \sum_{k=1}^n \| u_j\|_{H^1(A_k^n)}^2 \leq \|u_j\|_{H^1(B({p_j\over \e_j},n)^c)}^2 \leq 4\nu^2,
    $$
we can find a $k_j^n\in \{ 1,2,\cdots,n\}$ with the desired property.

\smallskip

\noindent
For $k_j^n$ given in Step 1, we set
    $$  v_j^n(x)=\widetilde\zeta_{n(k_j^n+1)}(|x-{p_j\over \e_j}|)u_j(x),
    $$
where $\widetilde\zeta_R(s)$ is given in \eqref{eq:2.4}.

\smallskip

\noindent
{\sl Step 2:} $\displaystyle \limsup_{j\to\infty} |J_\e'(v_j^n)v_j^n| \leq C_n$, 
where $C_n>0$ satisfies $C_n\to 0$ as $n\to\infty$.

\smallskip

\noindent
We compute
    \bea
    && (J_{\e_j}'(u_j)-J_{\e_j}'(v_j^n)) v_j^n \\
    &=& \int_{\R^N} (\nabla u_j-\nabla v_j^n)\nabla v_j^n+ V(\e_jx)(u_j-v_j^n)v_j^n\, dx 
        -\half (\calD'(u_j)-\calD'(v_j^n))v_j^n \\
    &=& (I)+(II).
    \eea
Clearly,
    $$  |(I)| \leq C\|u_j\|_{H^1(A_{k_j^n}^n)}^2 \leq C{4\nu^2\over n}.
    $$
By Corollary \ref{S:2.6} (iii), 
    $$
    |(II)| \leq C_1\|u_j\|_{H^1(A_{k_j^n}^n)}^{N+\alpha\over N} +C_2 m_n 
    \leq C_1 ({4\nu^2\over n})^{N+\alpha\over 2N} +C_2m_n.
    $$
Since $J_{\e_j}'(u_j)\to 0$, we have Step 2.

\smallskip

\noindent
{\sl Step 3:} $\displaystyle \limsup_{j\to \infty} \|u_j(x+{p_j\over \e_j})\|_{H^1(B(0,n(n+1))^c)} \leq C_n'$, 
where $C_n'>0$ satisfies $C_n'\to 0$ as $n\to\infty$.

\smallskip

\noindent
By Corollary \ref{S:2.3}, we have
    $$  \limsup_{j\to\infty} \|v_j^n\|_{H^1} \leq {1\over c_0}C_n.
    $$
Since $\| u_j\|_{H^1(B({p_j\over \e_j}, n(n+1))^c)} \leq \|v_j^n\|_{H^1}$, we have Step 3.

\smallskip

\noindent
{\sl Step 4:} $u_j(x+{p_j\over \e_j})\to u_0$ strongly in $H^1(\R^N)$.

\smallskip

\noindent
Set $w_j(x)=u_j(x+{p_j\over \e_j})$.  Since $w_j\wlimit u_0$ weakly in $H^1(\R^N)$, we have $w_j\to u_0$
strongly in $L^q_{loc}(\R^N)$ for any $q\in [2,{2N\over N-2})$.
By Step 3, we have
    $$  \lim_{n\to\infty}\limsup_{j\to\infty} \|w_j\|_{L^q(B(0,n(n+1))^c)} =0
    $$
for any $q\in [2,{2N\over N-2}]$.  Thus we have $w_j\to u_0$ strongly in $L^q(\R^N)$ for any $q\in [2,{2N\over N-2})$,
which implies
    $$  \calD'(w_j)\to \calD'(u_0) \quad \hbox{strongly in}\ H^{-1}(\R^N).
    $$
We compute
    \bea
    && \| \nabla w_j\|_2^2 + V(p_0) \| w_j\|_2^2 
    = \| \nabla w_j\|_2^2 + \int_{\R^N} V(\e_j x+p_j)w_j^2\, dx +o(1) \\
    &=& \| \nabla u_j\|_2^2 + \int_{\R^N} V(\e_j x)u_j^2\, dx +o(1) 
    = J_{\e_j}'(u_j)u_j + \half\calD'(u_j)u_j + o(1) \\
    &=& \half\calD'(u_j)u_j + o(1) 
    = \half\calD'(w_j)w_j + o(1) \\
    &=& \half\calD'(u_0)u_0.
    \eea
Similarly,
    \bea
    && \| \nabla u_0\|_2^2 + V(p_0) \| u_0\|_2^2 \\
    &=& \int_{\R^N} \nabla w_j\nabla u_0 + V(p_0)w_j u_0\, dx +o(1) \\
    &=& \int_{\R^N} \nabla w_j\nabla u_0 + V(\e_jx+p_0)w_j u_0\, dx +o(1) \\
    &=& \int_{\R^N} \nabla u_j\nabla u_0(x-{p_j\over\e_j}) + V(\e_j x)u_j u_0(x-{p_j\over\e_j})\, dx+o(1)\\
    &=& J_{\e_j}'(u_j)u_0(x-{p_j\over\e_j}) + \half\calD'(u_j)u_0(x-{p_j\over\e_j}) +o(1)\\
    &=& \half\calD'(u_j)u_0(x-{p_j\over\e_j}) + o(1) \\
    &=& \half\calD'(w_j)u_0+ o(1) \\
    &=& \half\calD'(u_0)u_0.
    \eea
Thus we have $\| \nabla w_j\|_2^2 + V(p_0) \| w_j\|_2^2 \to \| \nabla u_0\|_2^2 + V(p_0) \| u_0\|_2^2$, which
implies the conclusion of Step 4.

\smallskip

\noindent
{\sl Step 5:} Conclusion.

\smallskip

\noindent
We show that $p\in K$ and $u_0$ is a least energy solution of 
    \be \label{eq:4.3}
        -\Delta u+m_0 u=(I_\alpha*F(u))f(u)\quad \hbox{in}\  \R^N.
    \ee
By \eqref{eq:4.2}, we have for any $\varphi\in C_0^\infty(\R^N,\R)$
    $$  J_{\e_j}'(u_j)\varphi(x-{p_j\over \e_j}) \to 0,
    $$
from which we have
    $$  \int_{\R^N} \nabla u_0\nabla\varphi+V(p_0)u_0\varphi \, dx
        -\half \calD'(u_0)\varphi =0.
    $$
Thus $u_0$ is a solution of \eqref{eq:4.3}.  We also have
    $$  J_{\e_j}(u_j) \to \half \|\nabla u_0\|_2^2 +\half V(p_0)\|u_0\|_2^2 -\half\calD(u_0)
    $$
and we have $E(m_0)= \half \|\nabla u_0\|_2^2 +\half V(p_0)\|u_0\|_2^2 -\half\calD(u_0)$.  
Recalling $p_0\in \overline{K_{d_0}}$ and $\inf_{x\in \overline{K_{d_0}}} V(x)=m_0$, we have
$p_0\in K$ and $E(m_0)=L_{m_0}(u_0)$.  \QED

\bigskip

The following corollary gives an uniform estimate of $\norm{J'_\epsilon(u)}_{H^{-1}}$ in an
annular neighborhood of a set of expected solutions, which is one of the keys of our argument.

\medskip
\begin{corollary} \label{S:4.2}
Let $\nu_4 >0$ be given in Proposition \ref{S:4.1}.  Then for any $0 <\rho_1 < \rho_2 < \nu_4$ and for all $d \in (0, d_0)$ 
there exists $\delta_1>0$ such that for $\epsilon>0$ small
    $$  \norm{J_\epsilon'(u)}_{H^{-1}} \geq \delta_1
    $$
for all $u\in N_\epsilon(\nu_4)$ with $J_\epsilon (u)\in [E(m_0)- \delta_1, E(m_0) + \delta_1]$ and $(\wrho(u),\epsilon\beta(u))
    \in ([0,\rho_2]\times K_{d_0} ) \setminus([0,\rho_1]\times K_{d})$.
\end{corollary}

\claim Proof.
By contradiction, we assume that there exist 
$(\epsilon_j) \subset (0, 1]$ such that $\epsilon_j \to 0$, as $j \to + \infty$ and
$(u_j) \subset H^1(\R^N)$ such that 
    \begin{eqnarray*}
    && u_j \in  N_{\epsilon_j}(\nu_4),\\
    && (\wrho(u_j),\epsilon\beta(u_j)) \in ([0,\rho_2]\times K_{d_0} ) 
        \setminus([0,\rho_1]\times K_{d}).
    \end{eqnarray*}
and
    $$  J_{\epsilon_j}(u_j) \to E(m_0),  \quad 
        J'_{\epsilon_j}(u_j) \to 0 \quad \hbox{as} \ j \to + \infty.
    $$
 
By Proposition \ref{S:4.1},  there exist, up to a subsequence,   $(p_j) \subset \R^N$, $p \in K$  and $\omega_0 \in S_{m_0}$ such that $p_j \to p$ and
    $$
    \| u_j\bigl(\cdot + \frac{p_j}{\epsilon_j}\bigr) -\omega_0 \|_{H^1} \to 0
    \quad \hbox{as} \ j \to + \infty.
    $$
    Therefore we have $$
    \| u_j -\omega_0\bigl(\cdot - \frac{p_j}{\epsilon_j}\bigr) \|_{H^1} \to 0
    \quad \hbox{as} \ j \to + \infty
    $$
and thus $\epsilon_j \beta (u_j) \to p\in K$ and 
$\wrho(u_j) \to 0$. A contradiction follows. \QED

\bigskip

Finally we need to prove the following result.

\bigskip

\begin{proposition}\label{S:4.3}
For any $\epsilon\in (0,1]$ fixed, the Palais-Smale condition holds for $J_\epsilon$ in
the set  $N_\epsilon(\nu_4)$.
That is, if a sequence $(u_j) \subset N_\epsilon(\nu_4)$ satisfies
    \begin{eqnarray}
    &&J_\epsilon (u_j) \to c, \label{eq:4.4}\\
    &&J'_\epsilon (u_j) \to 0 \quad \hbox{as} \ \ j\to + \infty, \label{eq:4.5}
    \end{eqnarray}
for some constant $c\in\R$, 
then $(u_j)$ has a strong convergent subsequence in $H^1(\R^N)$.
\end{proposition}

\claim Proof.
Let $\epsilon >0$ be fixed and
a sequence $(u_j) \subset N_\epsilon(\nu_4)$ satisfy \eqref{eq:4.4}--\eqref{eq:4.5}.
Since $N_\epsilon (\nu_4)$ is bounded in $H^1(\R^N)$, we can assume that
$(u_j)$  weakly  converges in $H^1(\R^N)$
to some $u_0 \in H^1(\R^N)$,  up to subsequences.

As in Steps 1--3 of the proof of Proposition \ref{S:4.1}, we can show
    $$  \lim_{n\to\infty}\limsup_{j\to\infty} \| u_j\|_{H^1(B(\beta(u_j), n(n+1))^c)} = 0.
    $$
We note that $\beta(u_j)$ stays bounded as $j\to\infty$, since $\e\beta(u_j)\in K_{d_0}$ and 
$\e\in (0,1]$ is fixed.
Thus we have
    $$  \lim_{R\to\infty}\limsup_{j\to\infty} \| u_j\|_{H^1(B(0,R)^c)} = 0,
    $$
which implies $u_j\to u_0$ strongly in $H^1(\R^N)$.  \QED

\bigskip


\setcounter{equation}{0}

\section{\label{section:5} Deformation outside small balls}
We set
    \begin{equation}\label{eq:5.1}
    \nu_* = {1\over 6}\min\{ \nu_0, \nu_1, \cdots, \nu_4\}>0.
    \end{equation}
In this section, we construct a map $\tau_\e:\, N_\e(\nu_*)\to N_\e(5\nu_*)$ whose properties are given in the following
proposition.

\medskip

\begin{proposition}\label{S:5.1}
For sufficiently small $\e>0$, there exists a map \\
$\tau_\e:\, N_\e(\nu_*)\to N_\e(5\nu_*)$ such that
\begin{itemize}
\item[(i)] $\tau_\e:\, N_\e(\nu_*)\to N_\e(5\nu_*)$ is continuous.
\item[(ii)] $\tau_\e(u)=u$ if $u(x)=0$ for all $|x-\beta(u)|\geq {1\over\sqrt\e}$.
\item[(iii)] For all $u\in N_\e(\nu_*)$,
    \begin{eqnarray}
    &&\tau_\e(u)(x)=u(x) \quad \hbox{for}\ x\in B(\beta(u),{1\over\sqrt\e}), \label{eq:5.2}\\
    &&J_\e(\tau_\e(u)) \leq J_\e(u),                                         \label{eq:5.3}\\
    &&\|\tau_\e(u)\|_{H^1(B(\beta(u),{2\over\sqrt\e})^c)} \leq \widetilde \rho_\e, \label{eq:5.4}\\
    &&|\beta(\tau_\e(u))-\beta(u)| \leq 2R_0,                                \label{eq:5.5}\\
    &&\wrho(\tau_\e(u)) \leq 5\nu_*.                                            \label{eq:5.6}                 
    \end{eqnarray}
Here $\widetilde\rho_\e>0$ is independent of $u$ and satisfies $\widetilde\rho_\e\to 0$ as
$\e\to 0$.
\end{itemize}
\end{proposition}

\medskip

We note that this type of operators were introduced in Byeon and Tanaka \cite{BT1, BT2, BT3} for nonlinear
Schr\"odinger equations and a related problem through a minimizing problem outside small balls, which is related to the
following boundary value problem:
    \bea
    &&-\Delta v +V(\e x)v = f(v) \quad \hbox{in}\ B(\beta(u), {1\over \sqrt\e})^c,\\
    &&\ v=u \quad \hbox{on}\ \partial B(\beta(u), {1\over \sqrt\e}),\\
    &&\ \| v\|_{H^1(B(\beta(u), {1\over \sqrt\e})^c)} \leq 5\nu_*.
    \eea
That is, $\tau_\e(u)$ is given as a solution of the above boundary value problem. 
In this argument, to obtain unique solvability and continuity of the operator, strict convexity of the corresponding functional 
is important.   Because of non-local nonlinearity, it can be verified for our
problem just under the condition $\limsup_{t\to 0}|{f(t)\over t}| <\infty$, that is, $p\geq 2$ for $f(s)=|s|^{p-2}s$.

Here we take another approach to construct $\tau_\e$ and it will be constructed using a special deformation flow outside
small balls.

To make our construction of $\tau_\e$ clearly, we fix the center of balls and 
work in a fixed space $H^1(B(0,{1\over \sqrt\e})^c)$.  First we note
    $$  J_\e(u) = \wJ_\e(\e\beta(u), u(x+\beta(u))),
    $$
where 
    $$  \wJ_\e(p,v)=\half\|\nabla v\|_2^2 +\half\int_{\R^N} V(\e x+p)v^2\, dx -\half\calD(v).
    $$
We denote restriction of $u(x+\beta(u))$ on the set $\Bonec$ by $\Theta_\e(u)$.  We note
that for all $u\in N_\e(\nu_*)$
    $$  \Theta_\e(u) \in Y_\e\equiv \{ v\in H^1(\Bonec);\, \|v\|_{H^1(\Bonec)}< {3\over 2}\nu_*\}
        \quad \hbox{for $\e>0$ small}.
    $$
In what follows, we construct a vector field $V(v):\, Y_\e \to H_0^1(\Bonec)$ with special properties.
We define $K_{1\e}(v)$, $K_{2\e}(v):\, H^1(B(0,{1\over\sqrt\e})^c)\to \R$ by
    $$  K_{1\e}(v) = \|v\|_{H^1(\Bonec)}^2,\quad 
        K_{2\e}(v) = \|v\|_{H^1(\Btwoc)}^2,
    $$
which will play important roles to construct $\tau_\e$.
We also set \\
$M=\sup_{\omega\in S_{m_0}, s\in [1-s_0,1+s_0]} \| \omega(x/s)\|_{H^1}+1$.

\begin{lemma}\label{S:5.2}
For any $v\in Y_\e$, there exists a $\calV_v\in H_0^1(\Bonec)$ such that
\begin{itemize}
\item[(i)] 
    \be\label{eq:5.7}
    \|\calV_v\|_{H^1} \leq 2\nu_*.
    \ee
\item[(ii)] If $u\in H^1(\R^N)$ satisfies 
    \be\label{eq:5.8}
    \|u\|_{H^1}\leq M  \quad \hbox{and}\quad  u=v \ \hbox{in}\  \Bonec,
    \ee
then for any $p\in K_{2d_0}$
    \be\label{eq:5.9}
    \wJ_\e'(p,u)\calV_v \geq a(\|v\|_{H^1(\Btwoc)}^2-\rho_\e).
    \ee
Here $a>0$ is independent of $\e$, $p$, $u$, $v$ and $\rho_\e$ is independent of $p$, $u$, $v$ and
satisfies
    $$  \rho_\e\to 0 \quad \hbox{as}\ \e\to 0.
    $$
\item[(iii)] For any $v\in Y_\e$, 
    \be \label{eq:5.10}
    K_{1\e}'(v)\calV_v,\  K_{2\e}'(v)\calV_v\geq a(\|v\|_{H^1(\Btwoc)}^2 -\rho_\e).
    \ee
Here $a>0$ and $\rho_\e>0$ are as in (ii).
\end{itemize}
\end{lemma}

\claim Proof.
Proof is divided into 3 steps.

We denote by $n_\e$ the largest integer less than $1/\e^{1/4}$.  
\smallskip

\noindent
{\sl Step 1:} For any $v\in Y_\e$, there exists $k\in\{ 1,2,\cdots, n_\e\}$ such that
    $$  \|v\|_{H^1(|x|\in [{1\over\sqrt\e}+{k-1\over\e^{1/4}},{1\over\sqrt\e}+{k\over\e^{1/4}}])}^2
        < ({3\over 2})^2{\nu_*^2\over n_\e}.
    $$

\smallskip

\noindent
In fact, we have
    $$  \sum_{k=1}^{n_\e} \|v\|_{H^1(|x|\in [{1\over\sqrt\e}+{k-1\over\e^{1/4}},{1\over\sqrt\e}+{k\over\e^{1/4}}])}^2
        \leq \|v\|_{H^1(\Bonec)}^2 \leq ({3\over 2})^2\nu_*^2.
    $$
Thus Step 1 holds.

\smallskip

\noindent
For a $k$ given in Step 1, we set
    $$  \calV_v(x)=\widetilde \zeta_{{1\over\sqrt\e}+{k\over\e^{1/4}}}(|x|)v(x)\in H_0^1(\Bonec).
    $$
Here $\widetilde \zeta_R(s)$ is defined in \eqref{eq:2.4}.  
For $\e>0$ small, clearly we have the property \eqref{eq:5.7}.

\smallskip

\noindent
{\sl Step 2:} $\calV_v$ satisfies \eqref{eq:5.9}.

\smallskip

\noindent
We compute
    \begin{eqnarray*}
    &&\wJ_\e'(p,u)\calV_v \\
    &=& \int_{\R^N}\nabla v\nabla(\zR v)\, dx +\int_{\R^N}V(\e x+p)\zR v^2\, dx \\
    && -\half\calD'(u)(\zR v)\\
    &=& \int_{|x|>{1\over\sqrt\e}+{k\over \e^{1/4}}} |\nabla v|^2+V(\e x+p)v^2\,dx \\
    && + \int_{|x|\in [{1\over\sqrt\e}+{k-1\over \e^{1/4}},{1\over\sqrt\e}+{k\over \e^{1/4}}]} 
            \nabla v\nabla(\zR v) +V(\e x+p)\zR v^2\, dx \\
    && -\half\calD'(u)(\zR v)\\
    &\geq& c\| v\|_{H^1(|x|>{1\over\sqrt\e}+{k\over\e^{1/4}})}^2 +(I) -(II),
    \end{eqnarray*}
where $c=\min\{ 1,\underline V\}$.

We can easily see that
    $$  |(I)| \leq c'\|v\|_{H^1(|x|\in [{1\over\sqrt\e}+{k-1\over \e^{1/4}},{1\over\sqrt\e}+{k\over \e^{1/4}}])}^2
        \leq C'{\nu_*^2\over n_\e}.
    $$
By Corollary \ref{S:2.6} (iii),
    \begin{eqnarray*}
    &&|(II)-\calD'(\zR v)(\zR v)| \\
    \leq&& C_1\|v\|_{H^1(|x|\in [{1\over\sqrt\e}+{k-1\over \e^{1/4}},{1\over\sqrt\e}+{k\over \e^{1/4}}])}^{N+\alpha\over N} 
            +C_2 m_{1/\e^{1/4}}\\
    \leq&& C_1\left(({3\over 2})^2{\nu_*^2\over n_\e}\right)^{N+\alpha\over 2N} +C_2m_{1/\e^{1/4}}.
    \end{eqnarray*}
By Proposition \ref{S:2.1}, using $\sigma(t)$ defined in \eqref{eq:2.3}, we have
    \begin{eqnarray*}
    && |\calD'(\zR v)(\zR v)| \\
    &\leq& C \| F(\zR v)\|_{2N\over N+\alpha} \| f(\zR v)(\zR v)\|_{2N\over N+\alpha} \\
    &\leq& C' \sigma(\| \zR v\|_{H^1})^{N+\alpha\over N} \\
    &\leq& C'' \sigma(\| v\|_{H^1(|x|>{1\over\sqrt\e}+{k\over \e^{1/4}})})^{N+\alpha\over N}
     + C''' \sigma(\| v\|_{H^1(|x| \in [{1\over\sqrt\e}+{k-1\over \e^{1/4}},{1\over\sqrt\e}+{k\over \e^{1/4}}])})^{N+\alpha\over N} \\
    &\leq& C'' \sigma(\| v\|_{H^1(|x|>{1\over\sqrt\e}+{k\over \e^{1/4}})})^{N+\alpha\over N}
     + C''' \sigma\left(\sqrt{\nu_*^2\over n_\e}\right)^{N+\alpha\over N}.
    \end{eqnarray*}
Choosing $\nu_4>0$ smaller if necessary, we may assume that for some $a>0$ 
    $$  ct^2 - C''\sigma(t)^{N+\alpha\over N} \geq a t^2 \quad \hbox{for}\ t\in [0,{3\over 2}\nu_*].
    $$
Thus, we have
    \begin{eqnarray*}
    \widetilde J_\e'(p,u)\calV_v &\geq& a\|v\|_{H^1(|x|> {1\over\sqrt\e}+{k\over\e^{1/4}})}^2 -c_\e \\
        &\geq& a\|v\|_{H^1(|x|>{2\over\sqrt\e})}^2 -c_\e,
    \end{eqnarray*}
where $c_\e= C'{\nu_*^2\over n_\e}+C_1\left(({3\over 2})^2 {\nu_*^2\over n_\e}\right)^{N+\alpha\over N}
+ C''' \sigma\left(\sqrt{\nu_*^2\over n_\e}\right)^{N+\alpha\over N}+ C_2m_{1/\e^{1/4}}$. 
Setting $\rho_\e={1\over a}c_\e$, we have \eqref{eq:5.9}.

\smallskip

\noindent
{\sl Step 3:}  $\calV_v$ satisfies \eqref{eq:5.10}.

\smallskip

\noindent
Since
    \begin{eqnarray*}
    K_{1\e}'(v)\calV_v &=& 2\int_{\Bonec}\nabla v\nabla(\zR v)+\zR v^2\, dx,\\
    K_{2\e}'(v)\calV_v &=& 2\int_{\Btwoc} |\nabla v|^2  + v^2\, dx,
    \end{eqnarray*}
\eqref{eq:5.10} can be shown as in Step 2.  \QED

\medskip

\noindent
We define
    \bea
    \whY &=& \{ v\in Y_\e;\, K_{2\e}(v) > 3\rho_\e\} \\
    &=& \{ v\in H^1(\Bonec);\, \|v\|_{H^1(\Bonec)} < {3\over 2}\nu_*,\ \|v\|_{H^1(\Btwoc)}^2 > 3\rho_\e\}.
    \eea
Setting $\widetilde \calV_v=\calV_v/\|\calV_v\|_{H^1}$, we have

\begin{lemma}\label{S:5.3}
For any $v\in \whY$, there exists $\widetilde\calV_v\in H_0^1(\Bonec)$
such that 
\begin{itemize}
\item[(i)] $\|\widetilde\calV_v\|_{H^1}\leq 1$.
\item[(ii)] $\wJ_\e'(p,u)\widetilde\calV_v$, $K_{1\e}'(v)\widetilde\calV_v$, $\displaystyle K_{2\e}'(v)\widetilde\calV_v
>{a\rho_\e\over \nu_*}$ for $u\in H^1(\R^N)$ with \eqref{eq:5.8} and $p\in K_{2d_0}$.
\end{itemize}
\end{lemma}

Clearly for any $v\in \whY$, there exists an open neighborhood $U_v$ in 
$\whY$ such that
    $$  \wJ_\e'(p,\widehat u)\widetilde\calV_v, \ K_{1\e}'(\widehat v)\widetilde\calV_v, \ 
        K_{2\e}'(\widehat v)\widetilde\calV_v   >{a\rho_\e\over \nu_*}
    $$
hold for $\widehat v\in U_v$, where $\widehat u$ and $\widehat v$ satisfies \eqref{eq:5.8}
and $p\in K_{2d_0}$.

Using partition of unity, we have

\begin{proposition}\label{S:5.4}
There exists a locally Lipschitz vector field $V(v):\, \whY\to H_0^1(\Bonec)$ such that 
for all $v\in \whY$ 
\begin{itemize}
\item[(i)] $\|V(v)\|_{H^1}\leq 1$.
\item[(ii)] $\wJ_\e'(p,u)V(v)$, $K_{1\e}'(v)V(v)$, $\displaystyle K_{2\e}'(v)V(v)>{a\rho_\e\over \nu_*}$ 
for $u\in H^1(\R^N)$ with \eqref{eq:5.8} and $p\in K_{2d_0}$.
\end{itemize}
\end{proposition}

We choose a function $\varphi(r)\in C^\infty(\R,\R)$ such that
    \begin{eqnarray*}
    &&\varphi(r)=\begin{cases} 1 &\hbox{for $r\geq 4\rho_\e$}, \\ 0 &\hbox{for $r\leq 3\rho_\e$}, \end{cases}\\
    &&\varphi(r)\in [0,1] \ \hbox{for all}\ r\in \R.
    \end{eqnarray*}
We consider the following ODE in $Y_\e$: for $v\in Y_\e$
    \begin{equation}\label{eq:5.11}
    \left\{
    \begin{array}{l}
        \displaystyle {dw\over d\tau}= -\varphi(K_{2\e}(w)) V(w),\\
        w(0,v)=v.
    \end{array}
    \right.
    \end{equation}

\begin{lemma}\label{S:5.5}
Let $w(\tau,v)$ be the unique solution of \eqref{eq:5.11}.  Then
\begin{itemize}
\item[(i)] $w(\tau,v)\in Y_\e$ for all $\tau\geq 0$ and $v\in Y_\e$.
\item[(ii)] Let $v\in Y_\e$ and suppose that $u\in H^1(\R^N)$ satisfies \eqref{eq:5.8} and $p\in K_{2d_0}$.
Then 
    $$  {d\over d\tau}\wJ_\e(p, \widetilde w_u(\tau,v)) \leq 0,
    $$
where
    $$  \widetilde w_u(\tau,v)(x)=\begin{cases}
            u(x)        &\hbox{for}\ x\in B(0,{1\over\sqrt\e}), \\
            w(\tau,v)   &\hbox{for}\ x\in B(0,{1\over\sqrt\e})^c. \\
        \end{cases}
    $$
\item[(iii)] For $i=1,2$,
    \bea  
    &&{d\over d\tau}K_{i\e}(w(\tau,v)) \leq 0 \quad \hbox{for all}\ v\in Y_\e, \\
    &&{d\over d\tau}K_{i\e}(w(\tau,v)) \leq -{a\rho_\e\over \nu_*} \quad \hbox{if}\ K_{2\e}(w(\tau,v))\geq 4\rho_\e.
    \eea
\item[(iv)]  There exists $T_\e>0$ such that for all $v\in Y_\e$
    $$  K_{2\e}(w(T_\e,v)) \leq 4\rho_\e.
    $$
\end{itemize}
\end{lemma}

\claim Proof.
First we show (iii).  Noting that $\varphi(K_{2\e}(w))>0$ implies \\
$\| w\|_{H^1(B(0,{1\over \sqrt\e})^c)} > 3\rho_\e$,
that is, $w\in \whY$,  we have from Proposition \ref{S:5.4}
    \begin{eqnarray}
    {d\over d\tau}K_{i\e}(w(\tau,v)) &=& -\varphi(K_{2\e}(w)) K_{i\e}'(w) V(w) \nonumber\\
        &\leq& -{a\rho_\e\over \nu_*}\varphi(K_{2\e}(w)). \label{eq:5.12}  
    \end{eqnarray}
Thus, (iii) holds.
By (iii), 
$K_{1\e}(w(\tau,v))=\| w(\tau,u)\|_{H^1(B(0,{1\over\sqrt\e})^c)}^2$ is non-\\ increasing.  
Thus $w(\tau,v)$ exists and satisfies $w(\tau,v)\in Y_\e$ for all $\tau\in [0,\infty)$.  
That is, $w(\tau,v):\, [0,\infty)\times Y_\e\to Y_\e$ is well-defined.
Thus (i) follows.

(iv) follows from (iii) easily.  Noting $w(\tau,v)|_{|x|={1\over\sqrt\e}}$ is independent of $\tau$, we have
$\widetilde w_u(\tau,v)\in H^1(\R^N)$.  Therefore (ii) follows from Proposition \ref{S:5.4} (ii).  \QED

\medskip

Now we define for $u\in N_\e(\nu_*)$
    $$  W(\tau,u)(x) = \begin{cases}
            u(x)                                &\hbox{for}\ x\in B(\beta(u), {1\over \sqrt\e}).\\
            w(\tau,\Theta_\e(u))(x-\beta(u))    &\hbox{for}\ x\in B(\beta(u), {1\over \sqrt\e})^c
            \end{cases}
    $$
and
    \bea    
        \widehat W(\tau,u)(x) &=& W(\tau,u)(x+\beta(u)) \\
        &=&\begin{cases}
            u(x+\beta(u))                       &\hbox{for} \ x\in B(0,{1\over\sqrt\e}),\\
            w(\tau,\Theta_\e(u))(x) &\hbox{for} \ x\in B(0,{1\over \sqrt\e})^c.
        \end{cases}
    \eea
We have

\begin{lemma}\label{S:5.6}
For $\e>0$ small, we have
\begin{itemize}
\item[(i)] For any $\tau\in [0,\infty)$, 
    \begin{eqnarray}
    &&W(\tau,u)(x) = u(x) \quad \hbox{for}\ x\in B(\beta(u),{1\over\sqrt\e}), \label{eq:5.13}\\
    &&J_\e(W(\tau,u)) \leq J_\e(u), \label{eq:5.14}\\
    &&\|W(\tau,u)\|_{H^1(B(\beta(u),{1\over\sqrt\e})^c)} \leq 2\nu_*, \label{eq:5.15}\\
    && |\beta(W(\tau,u))-\beta(u)| \leq 2R_0, \label{eq:5.16}\\
    &&\wrho(W(\tau,u)) \leq 5\nu_*.         \label{eq:5.17}
    \end{eqnarray}
\item[(ii)] Moreover we have
    $$  \| W(T_\e,u) \|_{H^1(B(\beta(u),{2\over\sqrt\e})^c)} \leq 2\sqrt{\rho_\e}.
    $$
\end{itemize}
\end{lemma}

\claim Proof.
(i) Obviously we have \eqref{eq:5.13}.  \\
We have for $u\in N_\e(\nu_*)$
    \bea
    J_\e(W(\tau,u)) &=& \widetilde J_\e(\e\beta(u), W(\tau,u)(x+\beta(u))) \\
        &=& \widetilde J_\e(\e\beta(u), \widehat W(\tau,u)(x)).
    \eea
Thus we deduce \eqref{eq:5.14} from (ii) of Lemma \ref{S:5.5}.

For \eqref{eq:5.15}, we note that
    $$  \| u\|_{H^1(B(\beta(u),{1\over\sqrt\e})^c)} < 2\nu_* \quad \hbox{for}\ \e>0 \ \hbox{small}.
    $$
Since
    $$  \| W(\tau,u)\|_{H^1(B(\beta(u),{1\over\sqrt\e})^c)}^2 
        = \| \widehat W(\tau,u)\|_{H^1(B(0,{1\over\sqrt\e})^c)}^2 
        = K_{1\e}(w(\tau,\Theta_\e(u)))
    $$
is non-increasing by (iii) of Lemma \ref{S:5.5}, we have \eqref{eq:5.15}.

Next we show \eqref{eq:5.16} and \eqref{eq:5.17}.  We write $u\in N_\e(\nu_*)$ as $u(x)=\omega({x-p\over s}) + \varphi(x)$,
where $\omega\in S_{m_0}$, $p\in \R^N$, $s\in [1-s_0, 1+s_0]$ and $\|\varphi\|_{H^1} < \nu_*$.
By Lemma \ref{S:3.5}, we have
    \begin{equation}\label{eq:5.18}
    |\beta(u)-p| \leq R_0.
    \end{equation}
Since
    \bea
    &&\| W(\tau,u)-u\|_{H^1} = \| W(\tau,u)-u\|_{H^1(B(\beta(u),{1\over\sqrt\e})^c)} \\
    &\leq& \| w(\tau,\Theta_\e(u))\|_{H^1(B(0,{1\over\sqrt\e})^c)} + 
            \| u\|_{H^1(B(\beta(u),{1\over\sqrt\e})^c)} 
    \leq 4\nu_*,
    \eea
we have
    \begin{equation}\label{eq:5.19}
    \| W(\tau,u)-\omega({x-p\over s})\|_{H^1}
    \leq \| W(\tau,u)-u\|_{H^1} + \| u-\omega({x-p\over s})\|_{H^1} 
    \leq 5\nu_*.  
    \end{equation}
Thus by our choice of $\nu_*$, we have
    $$  |\beta(W(\tau,u)) -p| \leq R_0,
    $$
from which and \eqref{eq:5.18} we have \eqref{eq:5.16}.  Clearly \eqref{eq:5.17} follows from \eqref{eq:5.19}.

(ii) Since
    $$  \| W(\tau,u)\|_{H^1(B(\beta(u),{2\over\sqrt\e})^c)}^2
        = \| w(\tau,\Theta_\e(u))\|_{H^1(B(0,{2\over\sqrt\e})^c)}^2
        = K_{2\e}(w(\tau,\Theta_\e(u))),
    $$
(ii) follows from (iv) of Lemma \ref{S:5.5}.  \QED

\medskip

Now we can give a proof to Proposition \ref{S:5.1}.

\medskip

\claim Proof of Proposition \ref{S:5.1}.
We set for $u\in N_\e(\nu_*)$
    $$  \tau_\e(u)= W(T_\e,u).
    $$
The desired properties \eqref{eq:5.2}--\eqref{eq:5.6} with $\widetilde\rho_\e=2\sqrt{\rho_\e}$
follow from Lemma \ref{S:5.6}.  \QED


\setcounter{equation}{0}
\section{\label{section:6}An invariant set and a deformation flow}

This section is devoted to develop a deformation argument. To this aim we introduce a neighborhood 
$\calXed$ of the set of expected solutions, which is positively invariant under a pseudo-gradient flow.

\medskip
\subsection{A pseudo-gradient flow}
For $\nu_*>0$ defined in \eqref{eq:5.1}, we fix $\rho_1$, $\rho_2>0$ such that
$0<\rho_1<\rho_2<\nu_*$ and we choose $\delta_0$, $\delta_1>0$ by Lemma \ref{S:3.7} and Corollary \ref{S:4.2}
with $d=\half d_0$.  We also set 
    $$  \delta_* = \min\{ {\delta_1\over 4}(\rho_2-\rho_1), \delta_0\}>0.
    $$
With these choices of $\nu_*$, $\delta_*$ we have for $\e>0$ small

\smallskip

\begin{itemize}
\item[(i)] $\Lambda(u)$ and $\beta(u)$ are well-defined on $N_\e(6\nu_*)$ and satisfy
    \bea
    &&\Lambda(u) \in [1-2s_0,1+2s_0], \\
    &&\e\beta(u) \in K_{2d_0}
    \eea
for all $u\in N_\e(6\nu_*)$.
\item[(ii)] By Lemma \ref{S:2.2}, for all $u(x)=\omega({x-p\over s})+\varphi(x)$ with $\omega\in S_{m_0}$,
$s\in [1-s_0, 1+s_0]$, $\|\varphi\|_{H^1}<6\nu_*$
    $$  |\beta(u)-p| < R_0.
    $$
\item[(iii)] By Lemma \ref{S:3.7},
    \be\label{eq:6.1}
    J_\e(u) \geq E(m_0) +\delta_0
    \ee
if $u\in N_\e(6\nu_*)$ satisfies $\e\beta(u)\in K_{{5\over 4}d_0}\setminus K_{{3\over 4}d_0}$.
\item[(iv)] By Corollary \ref{S:4.2},
    \be\label{eq:6.2}
    \| J_\e'(u)\|_{H^{-1}} \geq \delta_1,
    \ee
if $u\in N_\e(6\nu_*)$ satisfies $J_\e(u)\in [E(m_0)-\delta_1, E(m_0)+\delta_1]$ and \\
$(\wrho(u),\e\beta(u)) \in ([0,\rho_2]\times K_{d_0})\setminus ([0,\rho_1]\times K_{\half d_0})$.
\end{itemize}

\smallskip

\noindent
We define
$$  \calXed=\{ u\in N_\epsilon (\nu_*) \ | \,
J_\epsilon(u)\leq E(m_0)+\delta_*-{\delta_1\over 2}(\wrho(u)-\rho_1)_+\}.
$$
We shall try to find critical points of $J_\epsilon$ in $\calXed$.  We note that

\begin{itemize}
\item[(a)] By \eqref{eq:6.1}, for $u\in N_\epsilon (6\nu_*)$,  $\epsilon \beta(u) \not\in K_{{3\over 4}d_0}$
implies
\be\label{eq:6.3}
J_\epsilon(u)  \geq E(m_0)+\delta_0 > E(m_0)+\delta_*,
\ee
which implies $u\not\in\calXed$.  Thus we have
    \be\label{eq:6.4}
    \e\beta(u) \in K_{{3\over 4}d_0} \quad \hbox{for all}\ u \in \calXed.
    \ee
\item[(b)]  For $u\in\calXed$, $\wrho(u) \geq \rho_2$ implies
    \begin{equation}\label{eq:6.5}
    J_\epsilon(u) \leq  E(m_0)+\delta_*-{\delta_1\over 2}(\rho_2-\rho_1) < E(m_0)-\delta_*.
    \end{equation}
\end{itemize}

\medskip

\noindent
The following deformation result plays an essential role to show the existence of 
critical points.

\medskip

\begin{proposition}\label{S:6.1}
    For any $c\in (E(m_0)-\delta_*,E(m_0)+\delta_*)$ and for any neighborhood $U$ of
    $\calK_c\equiv \{ u\in \calXed;\, \Je'(u)=0,\ \Je(u)=c\}$ ($U=\emptyset$ if
    $\calK_c=\emptyset$), there exist $r>0$ with $(c-r,c+r)\subset (E(m_0)-\delta_*,E(m_0)+\delta_*)$
and a deformation
    $\eta(\tau,u):\, [0,1]\times(\calXed\setminus U)\to\calXed$ such that
    \begin{itemize}
        \item[(i)] $\eta(0,u)=u$ for all $u$.
        \item[(ii)] $\eta(\tau,u)=u$ for all $\tau\in [0,1]$ if $\Je(u)\not\in [E(m_0)-\delta_*,
        E(m_0)+\delta_*]$.
        \item[(iii)] $\Je(\eta(\tau,u))$ is a non-increasing function of $\tau$ for all $u$.
        \item[(iv)] $\Je(\eta(1,u))\leq c-r$ for all $u\in \calXed\setminus U$ satisfying
        $\Je(u)\leq c+r$.
    \end{itemize}
\end{proposition}

\medskip

\claim Proof.
We consider a deformation flow defined by
\begin{equation}\label{eq:6.6}
\begin{cases}
{d\eta\over d\tau}=-\phi(\eta) {\calV(\eta)\over\norm{\calV(\eta)}_{H^1}},\\
\eta(0,u)=u,
\end{cases}
\end{equation}
where $\calV(u):\, \{ u\in H^1(\R^N);\, \Je'(u)\not=0\}\to H^1(\R^N)$ is a
locally Lipschitz continuous pseudo-gradient vector field satisfying
$$  \norm{\calV(u)}_{H^1}\leq \norm{\Je'(u)}_{H^{-1}},\quad
\Je'(u)\calV(u)\geq \half\norm{\Je'(u)}_{H^{-1}}^2
$$
and $\phi(u):\, H^1(\R^N)\to [0,1]$ is a locally Lipschitz
continuous function, which is defined in the standard way.
We require that $\phi(u) $ satisfies 
    $$  \phi(u)= \begin{cases}
        0 &\hbox{if}\  \Je(u)\not\in [E(m_0)-\delta_*,E(m_0)+\delta_*]. \\
        1 &\hbox{if}\  \Je(u)\in [c-\kappa, c+\kappa]\setminus {\cal N},
        \end{cases}
    $$
where $[c-\kappa, c+\kappa]\subset (E(m_0)-\delta_*,E(m_0)+\delta_*)$ and ${\cal N}$ is 
an open set such that $\calK_c\subset {\cal N}\subset U$.

We consider the flow defined by \eqref{eq:6.6}. 
The properties (i)-(iii) follows by standard arguments from the definition \eqref{eq:6.6} and since
$\phi(u)=0$ if $\Je(u)\not\in [E(m_0)-\delta_*,E(m_0)+\delta_*]$ . Clearly also since, by Proposition \ref{S:4.3},  $\Je$ satisfies the Palais-Smale condition for fixed $\epsilon >0$, property (iv) is standard. Thus to end the proof we just need to show that

\be\label{eq:6.7}
\eta(\tau,\calXed)\subset\calXed \quad \hbox{for all}\ \tau\geq 0,
\ee
namely that $\calXed $ is positively invariant under our flow. First note that because of  property  (iii), \eqref{eq:6.3} implies that for
$u\in\calXed$, $\eta(t)=\eta(t,u)$ remains in the set $\{u;\, \epsilon \beta(u)\in K_{{3\over 4}d_0}\}$. 
Also, because of property (ii) and \eqref{eq:6.5}, for  $u\in\calXed$, $\eta(t)$ remains in
$N_\epsilon(\nu_*)$. Thus to show \eqref{eq:6.7}, we just need to prove that the property
    $$  \Je(u)\leq E(m_0)+\delta_*-{\delta_1\over 2}(\wrho(u)-\rho_1)_+
    $$
is stable under the deformation. For this it suffices to show
that for a solution $\eta(\tau)$ of \eqref{eq:6.6},
if $0<s<t<1$ satisfies
\begin{eqnarray*}
    &\wrho(\eta(\tau)) \in [\rho_1,\rho_2] \quad \hbox{for all}\ \tau\in [s,t],\\
    &\Je(\eta(s))\leq E(m_0)+\delta_*-{\delta_1\over 2}(\wrho(\eta(s))-\rho_1),
\end{eqnarray*}
then
$$  \Je(\eta(t))\leq E(m_0)+\delta_* -{\delta_1\over 2}(\wrho(\eta(t))-\rho_1).
$$
We note that $(\wrho(\eta(\tau)), \epsilon\beta(\eta(\tau)))\in [\rho_1,\rho_2]\times K_{d_0}
\subset \bigl([0, \rho_2] \times K_{d_0} \bigr)
\setminus 
\bigl([0, \rho_1] \times K_{\half d_0}\bigr)$
for all $\tau\in[s,t]$.
Thus by Corollary \ref{S:4.2}, we have for $\tau\in [s,t]$
    \bea  
    {d\over d\tau}\Je(\eta(\tau)) &=& \Je'(\eta){d\eta\over d\tau}
    = -\phi(\eta)\Je'(\eta){\calV(\eta)\over\norm{\calV(\eta)}_{H^1}}\\
    &\leq& -\half \|J_\e'(\eta)\|_{H^{-1}} 
    \leq -\phi(\eta){\delta_1\over 2}
    \eea
and
\be\label{eq:6.8}
\Je(\eta(t)) \leq \Je(\eta(s)) -{\delta_1\over 2}\int_s^t \phi(\eta(\tau))\, d\tau.
\ee
On the other hand,
\be\label{eq:6.9}
\norm{\eta(t)-\eta(s)}_{H^1} \leq \int_s^t \norm{\frac{d \eta}{d \tau}}_{H^1} \, d\tau
\leq \int_s^t \phi(\eta(\tau))\, d\tau.
\ee
By \eqref{eq:6.8}--\eqref{eq:6.9}, and using the fact that $\abs{\wrho(\eta(t))-\wrho(\eta(s))}\leq
\norm{\eta(t)-\eta(s)}_{H^1}$, we have
\begin{eqnarray*}
    \Je(\eta(t)) &\leq& \Je(\eta(s)) -{\delta_1\over 2}\norm{\eta(t)-\eta(s)}_{H^1}\\
    &\leq& \Je(\eta(s)) -{\delta_1\over 2}\abs{\wrho(\eta(t))-\wrho(\eta(s))}\\
    &\leq& E(m_0)+\delta_*-{\delta_1\over 2}(\wrho(\eta(s))-\rho_1)
    -{\delta_1\over 2}\abs{\wrho(\eta(t))-\wrho(\eta(s))}\\
    &\leq& E(m_0)+\delta_*-{\delta_1\over 2}(\wrho(\eta(t))-\rho_1).
\end{eqnarray*}
Thus $\eqref{eq:6.7}$ holds and the proof of the proposition is completed.
\QED

\medskip

\noindent
Proposition \ref{S:6.1} enables us to estimate the multiplicity of critical points using the
relative category.


\medskip

\subsection{\label{subsection:6.2} Two maps $\Phi_\epsilon$ and $\Psi_\epsilon$}
For $a\in\R$, we set
$$  \calXed^a=\{ u\in\calXed \ | \  \Je(u)\leq a\}.
$$
For $\hdelta\in (0,\delta_*)$ small, using relative category, we shall estimate the change of topology between $\calXp$ and
$\calXm$.
We introduce two maps:
\begin{eqnarray*}
    &&\Phi_\epsilon:\, ([1-s_0,1+s_0]\times K, \{ 1\pm s_0\}\times K)
    \to (\calXp, \calXm),\\
    &&\Psi_\epsilon:\, (\calXp, \calXm) \to\\
    &&\qquad ([1-s_0,1+s_0]\times K_{d_0}, ([1-s_0,1+s_0]\setminus\{ 1\})\times  K_{d_0})
    \hss.
\end{eqnarray*}
Here we use notation from algebraic topology: $f:\, (A,B)\to (A',B')$ means $B\subset A$,
$B'\subset A'$, $f:\, A\to A'$ is continuous and $f(B)\subset B'$.

\medskip

\noindent
{\sl Definition of $\Phi_\epsilon$}: \\
We fix a positive least energy solution $\omega_0(x)\in S_{m_0}$ 
and set for $(s,p)\in [1-s_0,1+s_0]\times K$
    $$  \Phi_\e(s,p)=\zeta_{1\over 2\sqrt\e}(|x-{p\over \e}|)\omega_0({x-{p\over\e}\over s}).
    $$

\noindent
{\sl Definition of $\Psi_\epsilon$}:\\
We define $\Psi_\e$ as a composition of the following maps.
\begin{itemize}
\item[(1)] $N_\e(\nu_*)\to N_\e(6\nu_*);\, u\mapsto \zeta_{2\over\sqrt\e}(|x-\beta(u)|)\tau_\e(u)(x)$.
\item[(2)] $N_\e(6\nu_*)\to [1-s_0,1+s_0]\times K_{d_0}; \, v\mapsto (\psi(\Lambda(v)),\e\beta(v))$,
\end{itemize}
where $\tau_\e(u):\, N_\e(\nu_*)\to N_\e(5\nu_*)$ is given in Proposition \ref{S:5.1} and $\psi(s)\in C(\R,\R)$ is defined by
    $$  \psi(s)=\begin{cases}
            1-s_0   &\hbox{for}\ s<1-s_0, \\
            s       &\hbox{for}\ 1-s_0 \leq s\leq 1+s_0,\\
            1+s_0   &\hbox{for}\ s>1+s_0.
        \end{cases}
    $$
That is, 
    $$  \Psi_\e(u)=(\psi(\Lambda(\zeta_{2\over\sqrt\e}(|x-\beta(u)|)\tau_\e(u)(x))),\, \e\beta(\zeta_{2\over\sqrt\e}(|x-\beta(u)|)\tau_\e(u)(x))).
    $$
In what follows, we observe that $\Phi_\e$ and $\Psi_\e$ are well-defined for small $\hdelta>0$.

For $(s,p)\in [1-s_0,1+s_0]\times K$, we have as $\e\sim 0$
    \begin{eqnarray*}
    &&\e\beta(\zeta_{1\over 2\sqrt\e}(|x-{p\over\e}|) \omega_0({x-{p\over \e}\over s})) = p+o(1),\\
    &&J_\e(\zeta_{1\over 2\sqrt\e}(|x-{p\over\e}|) \omega_0({x-{p\over \e}\over s})) 
    = L_{m_0}(\omega_0({x-{p\over \e}\over s})) +o(1)\\
    &&\quad =\half g(s)\|\nabla \omega_0\|_2^2 + \half h(s)m_0\| \omega_0\|_2^2 +o(1)\\
    &&\quad =E(m_0) - \half (g(1)-g(s))\|\nabla \omega_0\|_2^2 - \half (h(1)-h(s))m_0\| \omega_0\|_2^2 +o(1).
    \end{eqnarray*}
Here $g(s)$, $h(s)$ are given in the proof of Lemma \ref{S:3.1}.  Choosing $\hdelta>0$ small so that
    $$  2\hdelta < \half(g(1)-g(1\pm s_0))\|\nabla \omega_0\|_2^2 + \half(h(1)-h(1\pm s_0))m_0\|\omega_0\|_2^2,
    $$
we find that $\Phi_\e$ is well-defined as a map $([1-s_0,1+s_0]\times K, \{ 1\pm s_0\}\times K)\to (\calXp,\calXm)$
for small $\e>0$.

Next we deal with $\Psi_\e$.  For $u\in \calXp$, we set $v_\e(x)=\zeta_{2\over \sqrt\e}(|x-\beta(u)|)\tau_\e(u)(x)$.
By (iii) of Proposition \ref{S:5.1},
    \begin{eqnarray*}
    J_\e(u) &\geq& J_\e(\tau_\e(u)) \\
        &=& J_\e(\tau_{2\over\sqrt\e}(|x-\beta(u)|)\tau_\e(u)) +o(1)\\
        &=& J_\e(v_\e) + o(1).
    \end{eqnarray*}
Here we use the fact that $\|\tau_\e(u)\|_{H^1(B(\beta(u),{2\over \sqrt\e})^c)} \leq \widetilde\rho_\e\to 0$ as
$\e\to 0$ uniformly in $u$.
Since $\supp v_\e \subset B(\beta(u),{2\over \sqrt\e}+2)$ and $\lim_{\e\to 0}\max_{x\in B(\beta(u), {2\over\sqrt\e}+2)}
|V(\e x)-V(\e\beta(u))|=0$, we have
    \begin{eqnarray*}
    J_\e(v_\e) &=& L_{V(\beta(u))}(v_\e) +o(1) \\
    &\geq& L_{m_0}(v_\e)+o(1) \\
    &\geq& E(m_0) (1-C_0(\Lambda(v_\e)-1)^2) + o(1).
    \end{eqnarray*}
Here we use Corollary \ref{S:3.2}.  Thus, for $u\in \calXm$,
    $$  E(m_0)(1-C_0(\Lambda(v_\e)-1)^2) \leq E(m_0) -\hdelta +o(1),
    $$
which implies 
    $$  \Lambda(v_\e)\not=1\quad \hbox{for}\  \e>0 \ \hbox{small}
    $$
and $\Psi_\e$ is well-defined as a map $(\calXp,\calXm)\to 
([1-s_0,1+s_0]\times K_{2d_0}, ([1-s_0,1+s_0]\setminus\{ 1\})\times K_{2d_0})$.
We also note that by \eqref{eq:6.4}, $\e\beta(v_\e)\in K_{{3\over 4}d_0}$.
Thus $\Psi_\e$ is well-defined as a map
    \bea
    &&(\calXp,\calXm) \\
    &&\to ([1-s_0,1+s_0]\times K_{{3\over 4}d_0}, ([1-s_0,1+s_0]\setminus\{ 1\})\times K_{{3\over 4}d_0}) \\
    &&\subset ([1-s_0,1+s_0]\times K_{d_0}, ([1-s_0,1+s_0]\setminus\{ 1\})\times K_{d_0}).
    \eea
\QED


\medskip

\noindent
The next proposition will be important to estimate the relative category $\cat(\calXp,\calXm)$.

\medskip

\begin{proposition}\label{S:6.2}
For $\e>0$ small,
    \begin{eqnarray*}
        &&\Psi_\epsilon\circ\Phi_\epsilon:\,([1-s_0,1+s_0]\times K, \{ 1\pm s_0\}\times K) \\
        &&\quad\to([1-s_0,1+s_0]\times K_{d_0},
        ([1-s_0,1+s_0]\setminus\{ 1\})\times K_{d_0})
    \end{eqnarray*}
    is homotopic to the embedding $j(s,p)=(s,p)$.  That is, there exists a continuous
    map
    $$  \eta:\, [0,1]\times[1-s_0,1+s_0]\times K \to [1-s_0,1+s_0]\times K_{d_0}
    $$
    such that
    \begin{eqnarray*}
        &&\eta(0,s,p)=(\Psi_\epsilon\circ\Phi_\epsilon)(s,p), \\
        &&\eta(1,s,p)=(s,p) \quad \hbox{for all}\ (s,p)\in [1-s_0,1+s_0]\times K,\\
        &&\eta(t,s,p)\in ([1-s_0,1+s_0]\setminus\{ 1\})\times K_d\\
        &&\qquad \hbox{for all}\ t\in [0,1] \ \hbox{and}\ (s,p)\in \{1\pm s_0\}\times K.
    \end{eqnarray*}
\end{proposition}

\claim Proof.
For $(s,p)\in [1-s_0,1+s_0]\times K$ and $\e>0$ small, we set
    $$  w_\e(x)=\Phi_\e(s,p)(x)=\zeta_{1\over 2\sqrt\e}(|x-{p\over\e}|)\omega({x-{p\over\e}\over s}).
    $$
We observe that 
$|\beta(w_\e)-{p\over \e}|\leq R_0$.  Thus $w_\e(x)=0$ in $B(\beta(w_\e),{1\over\sqrt\e})^c$ for $\e>0$ small,
which implies $\tau_\e(w_\e)=w_\e$.  Thus
    \begin{eqnarray*}
    (\Psi_\e\circ\Phi_\e)(s,p) &=& \Psi_\e(w_\e) = (\psi(\Lambda(w_\e)),\e\beta(w_\e)) \\
    &=& (s,p) +o(1),
    \end{eqnarray*}
which implies $\Psi_\e\circ\Phi_\e$ is homotopic to the embedding $j(s,p)=(s,p)$.  \QED

\bigskip

Differently from \cite{CL1, CL2},
(see \cite[Remark 4.3]{CJT}) 
we can not infer in general that 
\[
\cat(\calXp, \calXm) \geq  \cat(K, \partial K).
\]
Therefore in the work  it will be necessary to use the notions of category and cup-length for maps.


\medskip
\setcounter{equation}{0}
\section{\label{section:7} Proof of Theorem \ref{S:1.4}}

In order to prove our theorem, we shall need some topological tools that we now present for the reader convenience.
Following \cite{BW}, see also \cite{FW1, FW2}, we define

\medskip

\begin{definition} \label{S:7.1}
    Let $B \subset A$ and $B' \subset A'$ be topological spaces and $f:(A,B) \to (A',B')$ be a continuous map, that is $f : A \to A'$ is continuous and $f(B) \subset B'$. The category $\cat(f)$ of $f$ is the least integer $k \geq 0$ such that there exist open sets $A_0$, $A_1$, $\cdots$, $A_k$ with the following properties:
    \begin{itemize}
        \item[(a)] $A=A_0\cup A_1\cup\cdots\cup A_k$.
        \item[(b)] $B\subset A_0$ and there exists a map $h_0:\, [0,1]\times A_0\to A'$
        such that
        \begin{eqnarray*}
            &&h_0(0,x)=f(x) \qquad \hbox{for all}\ x\in A_0,\\
            &&h_0(1,x)\in B' \quad \qquad \hbox{for all}\ x\in A_0,\\
            &&h_0(t,x)=f(x) \qquad \hbox{for all}\ x\in B\ \hbox{and}\ t\in [0,1].
        \end{eqnarray*}
        \item[(c)] For $i=1,2,\cdots, k$, $f|_{A_i}:\, A_i\to A'$ is homotopic to a
        constant map.
    \end{itemize}
\end{definition}

We also introduce the cup-length of $f: (A,B) \to (A',B')$. Let $H^*$ denote Alexander-Spanier cohomology
with coefficients in the field $\F$.
We recall that the cup product $\smile$ turns $H^*(A)$ into a ring with unit $1_A$, and it turns $H^*(A,B)$
into a module over $H^*(A)$. A continuous map $f : (A,B) \to (A',B')$ induces a homomorphism
$f^* : H^*(A') \to H^*(A)$ of rings as well as a homomorphism $f^* : H^*(A',B') \to H^*(A,B)$ of
abelian groups.  We also use notation:
$$  \tilde H^n(A')=\begin{cases} 0 &\text{for $n=0$,}\\ H^n(A')&\text{for $n>0$.}\end{cases}
$$
For more details on algebraic topology we refer to \cite{Sp}.

\begin{definition} \label{S:7.2}
    For  $f:(A,B) \to (A',B')$ the cup-length, $\cuplength(f)$ is defined as follows;
    when $f^*:\, H^*(A',B')\to H^*(A,B)$ is not a trivial map, $\cuplength(f)$ is defined
    as the maximal integer $k\geq 0$ such
    that there exist elements $\alpha_1$, $\cdots$, $\alpha_k
    \in \tilde{H}^*(A')$  and $\beta\in H^*(A',B')$  with
    \begin{eqnarray*}  f^*(\alpha_1\smile\cdots\smile\alpha_k \smile \beta)
        &=& f^*(\alpha_1)\smile\cdots\smile f^*(\alpha_k)\smile f^*(\beta) \\
        &\not=& 0 \ \hbox{in}\ H^*(A,B).
    \end{eqnarray*}
    When $f^*=0:\, H^*(A',B')\to H^*(A,B)$, we define $\cuplength(f)=-1$.
\end{definition}

We note that $\cuplength(f)=0$ if $f^*\not=0:\, H^*(A',B')\to H^*(A,B)$ and $\tilde{H}^*(A')=0$.

As fundamental properties of $\cat(f)$ and $\cuplength(f)$, we have

\begin{proposition}\label{S:7.3}
    
    \begin{itemize}
        \item[(i)] For $f:\, (A,B)\to (A',B')$, $\cat(f)\geq \cuplength(f)+1$.
        \item[(ii)] For $f:\, (A,B)\to (A',B')$, $f':\, (A',B')\to (A'',B'')$,
        $$  \cuplength(f'\circ f)\leq \min\{\cuplength(f'), \cuplength(f)\}.
        $$
        \item[(iii)] If $f, g : (A,B) \to (A',B')$ are homotopic, then $\cuplength(f) = \cuplength(g).$
    \end{itemize}
\end{proposition}

\medskip

The proof of these statements can be found in \cite[Lemmas 2.6, 2.7]{BW}.
Finally we recall

\begin{definition} \label{S:7.4}
    For a set $(A,B)$, we define the relative category $\cat(A,B)$ and the relative
    cup-length $\cuplength(A,B)$ by
    \begin{eqnarray*}
        &&\cat(A,B)=\cat(id_{(A,B)}:\, (A,B)\to (A,B)),\\
        &&\cuplength(A,B)=\cuplength(id_{(A,B)}:\, (A,B)\to (A,B)).
    \end{eqnarray*}
    We also set
    $$  \cat(A)=\cat(A,\emptyset),\quad \cuplength(A)=\cuplength(A,\emptyset).
    $$
\end{definition}

\medskip

We also recall  the following topological lemma due to T. Bartsch (cf. \cite[Lemma 5.5]{CJT}) where we make use of the continuity property of Alexander-Spanier cohomology.

\medskip

\begin{lemma} \label{S:7.5}
    Let $K\subset \R^N$ be a compact set.  For a $d$-neighborhood $K_d=\{ x\in \R^N;\,
    \dist(x,K)\leq d\}$ and $I=[0,1]$, $\partial I=\{0,1\}$, we consider the inclusion
    $$  j:\, (I\times K,\partial I\times K)\to (I\times K_d,\partial I\times K_d)
    $$
    defined by $j(s,x)=(s,x)$.  Then for $d>0$ small,
    $$  \cuplength(j) \geq \cuplength(K).
    $$
\end{lemma}

Now we have all the ingredients to give the
\vspace{2mm}

\claim Proof of Theorem \ref{S:1.4}.
We observe that for $\epsilon>0$ small
\be\label{eq:7.1}
\#\{u\in \calXp\setminus\calXm;\, \Je'(u)=0\}
\geq \cat(\calXp,\calXm).
\ee
Using Proposition \ref{S:6.1}, \eqref{eq:7.1} can be proved in a standard way (c.f. \cite[Theorem 4.2]{FW1}).
\vspace{2mm}

By (i) of Proposition \ref{S:7.3}, we have
\be\label{eq:7.2}
\cat(\calXp,\calXm) \geq \cuplength(\calXp,\calXm)+1.
\ee
Since $\Psi_\epsilon\circ\Phi_\epsilon=\Psi_\epsilon\circ id_{(\calXp,\calXm)}\circ\Phi_\epsilon$,
it follows from (ii) of Proposition \ref{S:7.3} that
\begin{eqnarray}
\cuplength(\Psi_\epsilon\circ\Phi_\epsilon)
&\leq& \cuplength(id_{(\calXp,\calXm)}) \nonumber\\
&=& \cuplength(\calXp,\calXm). \label{eq:7.3}
\end{eqnarray}
By Proposition  \ref{S:6.2}, that $\Psi_\epsilon\circ\Phi_\epsilon$ is homotopic to the inclusion
    $$ j:\, (I\times K,\partial I\times K)\to (I\times K_{d_0},\partial I\times K_{d_0}). 
    $$
By (iii) of Proposition  \ref{S:7.3},
\be\label{eq:7.4}
\cuplength(\Psi_\epsilon\circ\Phi_\epsilon) = \cuplength(j).
\ee
At this point using \eqref{eq:7.2}--\eqref{eq:7.4} and 
recalling our choice of $d_0$, that is, Lemma \ref{S:7.5} holds for $d=d_0$
we deduce that
$$  \cat(\calXp,\calXm) \geq \cuplength(K)+1.
$$
Thus by \eqref{eq:7.1}, $\Je$ has at least $\cuplength(K)+1$ critical points in $\calXp\setminus\calXm$.
Recalling Proposition \ref{S:4.1}, this completes the proof of the Theorem.  \QED

\medskip

\claim Proof of Theorem \ref{S:1.1}.
Since $\cuplength(K)\geq 0$ for $K\not=\emptyset$, Theorem \ref{S:1.1} can be regarded as a special case of Theorem \ref{S:1.4}. 
\QED

\medskip

\claim Proof of Remark \ref{S:1.6}.
From Proposition \ref{S:4.1} we know that the critical points $u_{\e}^i$, $i= 1, \dots , \cuplength(K)+1$ satisfy
$$||u_{\e}^i (x) - \omega^i (x- x_\e^i) ||_{H^1} \to 0$$
where $\e x_\e^i = \e \beta(u_{\e}^i) + o(1) \to x_0^i \in K$ and $\omega^i \in S_{m_0}$. 
Thus $w_{\e}^i(x) = u_{\e}^i(x + x_\e^i)$ converges to $\omega^i \in S_{m_0}$.  
Each $w_{\e}^i$ converges to a least energy solution of \eqref{eq:1.8}.
\QED

\medskip


\setcounter{equation}{0}
\section{\label{section:8} Another proof of Theorem \ref{S:1.1}}
We can give another proof of Theorem \ref{S:1.1} via localized mountain pass argument.  We just give an outline of a proof.

\medskip

\claim Outline of a proof of Theorem \ref{S:1.1} via localized mountain pass. \\
Let $s_0>0$, $\nu_*>0$ be as in previous sections.  We may assume that $0\in K$.

It suffices to show for any $\delta\in (0, \nu_*)$ small that $J_\e(u)$ has a critical point $u\in N_\e(\delta)$
with $J_\e(u)\in [E(m_0)-\delta,E(m_0)+\delta]$ for $\e>0$ small.

Choosing $\delta>0$ smaller if necessary, as in Corollary \ref{S:4.2} we can show for some $\kappa_0>0$ independent of $\e$  that
for $\e>0$ small
    \begin{equation}\label{eq:8.1}
    \| J_\e'(u)\|_{H^{-1}} \geq \kappa_0 
    \end{equation}
for $u\in N_\e(\nu_*)$ with $(\wrho(u), \e\beta(u))\in ([0,\delta]\times K_{d_0})\setminus ([0,\half \delta]\times K_{{3\over 4}d_0})$
and $J_\e(u)\in [E(m_0)-\delta,E(m_0)+\delta]$.

If $J_\e(u)$ does not have any critical point with $J_\e(u)\in [E(m_0)-\delta,E(m_0)+\delta]$ in 
$N_\e(\delta)$, then, by Proposition \ref{S:4.3}, there exists $\kappa_\e\in (0,\kappa_0]$ such that
    \begin{equation}\label{eq:8.2}
    \| J_\e'(u)\|_{H^{-1}} \geq \kappa_\e \quad \hbox{for}\ u\in N_\e(\delta) \ \hbox{with}\  J_\e(u)\in [E(m_0)-\delta,E(m_0)+\delta].
    \end{equation}
Now we choose and fix $\omega_0\in S_{m_0}$.  We set
    $$  \gamma_\e(s)(x)=\zeta_{1\over 2\sqrt\e}(|x|) \omega_0({x\over s}):\, [1-s_0,1+s_0]\to H^1(\R^N).
    $$
We can easily see that
    \begin{eqnarray}
    &&\gamma_\e(s) \in N_\e({\delta\over 2}) \quad \hbox{for} \ s\in [1-s_0,1+s_0], \label{eq:8.3}\\
    &&\max_{s\in [1-s_0,1+s_0]} J_\e(\gamma_\e(s)) \to E(m_0) \quad \hbox{as}\ \e\to 0. \label{eq:8.4}
    \end{eqnarray}
If $J_\e(u)$ does not have any critical point with $J_\e(u)\in [E(m_0)-\delta, E(m_0)+\delta]$ in $N_\e(\delta)$, by
\eqref{eq:8.1}--\eqref{eq:8.2}, we can find a deformation flow $\eta:\, [0,T_\e]\times N_\e({\delta\over 2}) \to N_\e(\delta)$ such that
    \begin{eqnarray}
    &&\eta(T_\e,\gamma_\e(1\pm s_0)) = \gamma_\e(1\pm s_0), \label{eq:8.5} \\
    &&J_\e(\eta(T_\e,u)) \leq E(m_0) - {\widehat \delta\over 2} \ \hbox{for}\ u\in N_\e({\delta\over 2}) 
        \ \hbox{with}\ J_\e(u) \leq E(m_0)+ {\widehat \delta\over 2},      \label{eq:8.6}
    \end{eqnarray}
where $\widehat \delta = \half\min\{ \delta_0, \delta, {\kappa_0\delta\over 4}, E(m_0)-L_{m_0}(\omega({x\over 1\pm s_0}))\}>0$.
Here $\delta_0>0$ is given in Lemma \ref{S:3.7}.

In fact, the flow $\eta(\tau,u)$ is defined as a solution of
    \bea
    &&{d\eta\over d\tau}=- \varphi(\eta){V_\e(\eta)\over \| V_\e(\eta)\|_{H^1}}, \\
    &&\eta(0,u)=u,
    \eea
where $V_\e(u)$ is a pseudo-gradient vector field associated to $J_\e'(u)$ and $\varphi(u):\, N_\e(\delta)\to [0,1]$ is a locally
Lipschitz continuous function. 
With a suitable choice of $\varphi(u)$, we have 
    $$  \eta(\tau,u)=u \quad \hbox{if}\ J_\e(u)\not\in [E(m_0)-2\widehat\delta,E(m_0)+2\widehat\delta],
    $$
which implies \eqref{eq:8.5}.  Moreover, we have as long as 
$\eta(t,u)\in N_\e(\delta)$ satisfies $J_\e(\eta(\tau,u))\in [E(m_0)-\widehat\delta, E(m_0)+\widehat\delta]$, 
    \begin{eqnarray}
    &&\| {d\over d\tau}\eta(\tau,u)\|_{H^1}\leq 1,      \nonumber \\
    &&{d\over d\tau}J_\e(\eta(\tau,u)) \leq -\half \kappa_\e, \nonumber \\
    &&{d\over d\tau}J_\e(\eta(\tau,u)) \leq -\half \kappa_0 \quad \hbox{if}\ \eta(\tau,u)\in N_\e(\delta)\setminus N_\e(\half\delta).
                                                                                            \label{eq:8.7}
    \end{eqnarray}
We also note that by Lemma \ref{S:3.7}
    \begin{equation}\label{eq:8.8}
    J_\e(u) \geq E(m_0) +\delta_0 \geq E(m_0)+\widehat \delta 
    \ \hbox{if $u\in N_\e(\delta)$ satisfies $\e\beta(u)\in K_{d_0}\setminus K_{{3\over 4}d_0}$.}
    \end{equation}
For $u\in N_\e({\delta\over 2})$ with $J_\e(u)\in [E(m_0)-\widehat\delta,E(m_0)+\widehat\delta]$, we consider the flow
$\eta(\tau)=\eta(\tau,u)$.  As long as $\eta(\tau)\in N_\e(\delta)$ and 
$J_\e(\eta(\tau)) \in  [E(m_0)-\widehat\delta,E(m_0)+\widehat\delta]$
hold, we have ${d\over d\tau}J_\e(\eta(\tau)) \leq -\half \kappa_\e$, in particular, $J_\e(\eta(\tau))$ decreases at a certain rate.
By \eqref{eq:8.8}, 
if $\eta(\tau)$ leaves $N_\e(\delta)=\{ u\in \whS(\nu_*);\, \wrho(u) < \delta,\, \e\beta(u)\in K_{d_0}\}$, $\eta(\tau)$
gets out through the set $\{ u\in \whS(\nu_*);\, \wrho(u)=\delta\}$.  
Therefore, by \eqref{eq:8.7}, we deduce \eqref{eq:8.6} for large $T_\e$.

Thus, for sufficiently small $\e>0$, it follows from \eqref{eq:8.4} that \\
$\max_{s\in [1-s_0,1+s_0]} J_\e(\gamma_\e(s)) \leq E(m_0)+{\widehat\delta\over 2}$.  Thus $\widetilde \gamma_\e(s)=\eta(T_\e,\gamma_\e(s))$
satisfies
    \bea
    &&\widetilde\gamma_\e(1\pm s_0) =\gamma_\e(1\pm s_0), \\
    &&\max_{s\in [1-s_0,1+s_0]} J_\e(\widetilde \gamma_\e(s)) \leq E(m_0) -{\widehat \delta\over 2}.
    \eea
Let $\tau_\e(u):\, N_\e(\delta)\subset N_\e(\nu_*)\to N_\e(5\nu_*)$ be a map obtained in Proposition \ref{S:5.1}.  We set
    $$  \widehat \gamma_\e(s)= \tau_\e(\widetilde\gamma_\e(s)).
    $$
Thus we have
    \begin{eqnarray}
    &&J_\e(\widehat \gamma_\e(s)) \leq J_\e(\widetilde \gamma_\e(s)) \leq E(m_0)-{\widehat\delta\over 2}, \label{eq:8.9}\\
    && \| \widehat\gamma_\e(s)\|_{H^1(B(\beta(\widetilde\gamma_\e(s)),{2\over\sqrt\e})^c)} \leq\widetilde\rho_\e, \label{eq:8.10}\\
    &&\widehat \gamma_\e(1\pm s_0)=\widetilde\gamma_\e(1\pm s_0)=\gamma_\e(1\pm s_0). \label{eq:8.11}
    \end{eqnarray}
Here we use the fact that $\widetilde\gamma(1\pm s_0)|_{|x|\geq{1\over\sqrt\e}}=\gamma(1\pm s_0)|_{|x|\geq{1\over\sqrt\e}}=0$,
which implies \eqref{eq:8.11}.

Thus $\widehat \gamma_\e(s)$ satisfies $\beta(\widehat\gamma_\e(s))\in K_{d_0}$ for all $s\in [1-s_0,1+s_0]$ and it follows from
\eqref{eq:8.9}--\eqref{eq:8.10} that
    \begin{equation}\label{eq:8.12}
    L_{m_0}(\widehat \gamma_\e(s)) \leq J_\e(\widehat\gamma_\e(s)) + o(1) \leq E(m_0) -{\widehat\delta\over 2} +o(1)
    \end{equation}
uniformly in $s\in [1-s_0,1+s_0]$ as $\e\to 0$.

On the other hand, by \eqref{eq:8.11},
    $$  \Lambda(\widehat\gamma_\e(1\pm s_0)) = \Lambda(\gamma_\e(1\pm s_0))  = 1\pm s_0 + o(1) \quad \hbox{as} \ \e\to 0.
    $$
Thus for $\e>0$ small there exists $s_\e\in (1-s_0,1+s_0)$ such that 
    \begin{equation}\label{eq:8.13}
    \Lambda(\widehat\gamma_\e(s_\e))=1.
    \end{equation}
\eqref{eq:8.12} and \eqref{eq:8.13} contradicts with \eqref{eq:3.7}.  Thus $J_\e(u)$ has a critical point $u\in N_\e(\delta)$ satisfying
$J_\e(u)\in [E(m_0)-\delta,E(m_0)+\delta]$.  \QED

\bigskip

\noindent
{\bf 
Acknowledgments.}
Authors started this research during the second author's visit to 
Dipartimento di Meccanica, Matematica e Management, Politecnico di Bari in 2016
and the first author's visit to Department of Mathematics, Waseda University in 2017.
They would like to thank Politecnico di Bari and Waseda University for their kind 
hospitality.   
The first author is partially supported by 
INdAM-GNAMPA Project 2017  \lq\lq Metodi matematici per lo studio di fenomeni fisici nonlineari''.  
The second author is partially supported by JSPS Grants-in-Aid for Scientific 
Research (B) (25287025) and (B) (17H02855).


\end{document}